\documentclass[reqno]{amsart}
\usepackage{times}
\usepackage{graphics,epsfig,graphicx,color}
\usepackage{amsmath}
\usepackage{mathrsfs,amsfonts}
\usepackage{fullpage}
\usepackage{url}
\usepackage[debug,
letterpaper=true,
colorlinks=true,
linkcolor=red,
filecolor=green,
citecolor=red,
pdfpagemode=None]
{hyperref}

\graphicspath{{./FIGS/}}

\newtheorem{alg}{Algorithm}[section]
\newtheorem{rem}{Remark}[section]
\DeclareMathAlphabet{\itbf}{OML}{cmm}{b}{it}
 \DeclareMathAlphabet\mathbfcal{OMS}{cmsy}{b}{n}

\renewcommand{\hat}{\widehat}
\renewcommand{\tilde}{\widetilde}

\def\bx{{{\itbf x}}}
\def\bcu{{\boldsymbol{\it{U}}}}
\def\bu{{{\itbf U}}}
\def\bv{{\itbf V}}
\def\bw{{\itbf W}}
\def\cu{{\it{U}}}
\def\cp{{\it{P}}} 
\def\cb{{\it{b}}}
\def\bp{{\itbf P}}
\def\bb{{\itbf b}}
\def\bC{{\itbf C}}
\def\bn{{\boldsymbol{\nu}}}
\def\bD{{\itbf D}}
\def\bE{{\itbf E}}
\def\bI{{\itbf I}}
\def\bM{{\itbf M}}
\def\bL{{\itbf L}}
\def\cbL{\boldsymbol{\mathcal{L}}}
\def\tbL{\tilde{\bL}}
\def\ep{\varepsilon}
\def\la{\lambda}
\def\om{\omega}

\def\msP{{\mathscr{P}}}
\def\bbP{\pmb{\msP}}
\def\bA{{\itbf A}}

\def\tbA{{\pmb{\tilde{\mathcal{A}}}}}
\def\bbU{\pmb{\mathscr{U}}}
\def\bU{\pmb{{U}}}
\def\bP{\pmb{{P}}}
\def\bV{{\itbf V}}
\def\bW{{\itbf W}}
\def\bX{{\itbf X}}
\def\bR{{\itbf R}}
\def\bQ{{\itbf Q}}
\def\bM{{\itbf M}}

\def\bgamma{{\boldsymbol \gamma}}
\def\bgammahat{{\hat{\boldsymbol \gamma}}}
\def\balpha{{\boldsymbol \alpha}}
\def\bbeta{{\boldsymbol \beta}}

\begin{document}
\title{Untangling the nonlinearity in  inverse scattering with data-driven reduced order models}
\hyphenation{mamonov}
\hyphenation{mzaslavsky}

\author{Liliana Borcea}
\address{Department of Mathematics, University of Michigan,
Ann Arbor, MI 48109}
\email{borcea@umich.edu}

\author{Vladimir Druskin}
\address{Schlumberger-Doll Research Center, 1 Hampshire St.,
  Cambridge, MA 02139-1578}
\email{druskin1@slb.com}

\author{Alexander V. Mamonov}
\address{Department of Mathematics, University of Houston,
  Houston, TX 77204}
\email{mamonov@math.uh.edu}

\author{Mikhail Zaslavsky}
\address{Schlumberger-Doll Research Center, 1 Hampshire St.,
  Cambridge, MA 02139-1578}
\email{mzaslavsky@slb.com}

\begin{abstract}
The motivation of this work is  an inverse problem for the acoustic wave equation, where
an array of sensors probes an unknown medium with pulses and measures
the scattered waves. The goal of the inversion is to determine from these measurements
the structure of the scattering medium, modeled by a spatially varying
acoustic impedance function. Many inversion algorithms assume that the
mapping from the unknown impedance to the scattered waves is
approximately linear. The linearization, known as the Born
approximation, is not accurate in strongly scattering media, where the
waves undergo multiple reflections before they reach the sensors in
the array. Thus, the reconstructions of the impedance have numerous
artifacts. The main result of the paper is a novel, linear-algebraic algorithm that uses a 
reduced order model (ROM)  to map the data to those corresponding to the single scattering (Born) model. 
The ROM  construction is based only on the measurements at the sensors in the array. 
The ROM is a proxy for the wave propagator operator, 
that propagates the wave in the unknown medium over the duration of the time sampling interval. 
The output of the algorithm can be  input into any off-the-shelf inversion software that  
incorporates state of the art linear inversion algorithms to reconstruct the unknown 
acoustic impedance.
\end{abstract}
\maketitle



\thispagestyle{plain}

\section{Introduction}
\label{sec:intro}
Let us formulate the problem in a general setting, for  a hyperbolic
system of equations of the form 
\begin{align}
\partial_t \begin{pmatrix} \cp(t,\bx) \\ \bcu(t,\bx) \end{pmatrix} = 
\begin{pmatrix} 0 & - L_q \\ L_q^T & 0 \end{pmatrix}  
\begin{pmatrix} \cp(t,\bx) \\ \bcu(t,\bx) \end{pmatrix}, \quad \bx 
\in \Omega, ~ ~ t > 0,
\label{eq:in.1}
\end{align}
satisfied by a wave field with components $\cp(t,\bx)$
and $\bcu(t,\bx)$, in a simply connected domain
$\Omega \subset \mathbb{R}^d$ with piecewise smooth boundary $\partial
\Omega = \partial \Omega_a \cup \partial \Omega_i\,$ given by the
union of two sets: The first set is the accessible boundary $\partial
\Omega_a$,  where the measurements are made, and the second set is the inaccessible boundary 
$\partial \Omega_i$. In this paper we consider sound waves, where $\cp(t,\bx) \in \mathbb{R}$ corresponds 
to the acoustic pressure field and $\bcu(t,\bx) \in \mathbb{R}^d$ to
 the velocity field,  satisfying the boundary conditions
\begin{equation}
 \label{eq:in.3}
\bcu(t,\bx) \cdot \bn(\bx) = 0, \quad \bx \in \partial \Omega_a \quad \mbox{and} \quad  \cp(t,\bx) = 0, \quad \bx \in \partial \Omega_i,
\end{equation}
where $\bn(\bx)$ is the outer unit normal at $\partial \Omega_a$. However, the results can be extended  to 
other boundary conditions and to other waves that satisfy a system of form \eqref{eq:in.1}, 
such as electromagnetic and elastic waves. 

Note that  $\partial \Omega_i$  may be a true boundary or a fictitious one, for a truncation of
an infinite medium, in which case we can use causality to set the
condition \eqref{eq:in.3} at $\partial \Omega_i$, without affecting
the wave measured at the sensors on $\partial \Omega_a$, for the duration $t \le
t_{max}$.

The wave evolves in time $t$ starting from 
\begin{equation}
 \cp(0,\bx) = \cb(\bx), ~ ~  \bcu(0,\bx) = 0,
\label{eq:in.2}
\end{equation}
as described by the skew symmetric operator in \eqref{eq:in.1}, with
$L_q$ a first order partial differential operator in the $\bx$
variable, and $L_q^T$ its adjoint. The coefficients in these operators
depend linearly on a function $q(\bx)$, which is the unknown in the
inverse problem, to be determined from measurements of the wave,
``the data''. These are modeled by  continuously
differentiable measurement functions $\bM_j$ of the self-adjoint
operator $L_q L_q^T$,
\begin{equation}
\bD_j = \bM_j(L_q L_q^T), \quad j = 0, 1, \ldots, 2 n-1,
\label{eq:in5}
\end{equation}
with $\bD_j$ scalar or matrix valued, depending on the dimension $d$.

In the context of inverse
scattering for sound waves considered in this paper,  the data are gathered by  
a collection  (aka array) of sensors at $\partial \Omega_a$. These sensors act as both sources and receivers
that probe the medium with pulses and
measure the reflected pressure field. The acoustic system of wave equations takes
the form \eqref{eq:in.1} after a Liouville transformation of the
pressure and velocity fields, and the array measurements can be
written in the form \eqref{eq:in5}, as described in sections
\ref{sect:1D} and \ref{sect:MD}. If the medium has constant density,
there is a single unknown in the inverse problem, the wave speed $c(\bx)$. In
general, the medium has variable density $\rho(\bx)$, so we have two unknowns: the
wave speed $c(\bx)$ and the acoustic impedance $\sigma(\bx) = \rho(\bx) c(\bx)$.  

We
assume henceforth that the waves propagate through a medium with known
wave speed\footnote{In applications, we can only know the smooth wave
  speed in the medium that contains scattering inhomogeneities.  If
  these inhomogeneities have constant density, their impedance is a
  constant multiple of the wave speed, so by finding $\sigma(\bx)$ we
  can also determine the variations of $c(\bx)$.  Even if the perturbations
of $c(\bx)$ cannot be determined, their effect is mainly manifested in small travel time coordinate 
deformations, and our approach still suppresses multiple scattering artifacts, 
as illustrated in  section \ref{sect:num2D}.}, and the
unknown $q(\bx)$ in the inverse problem is the logarithm of the
acoustic impedance $\sigma(\bx)$.
This formulation of the inverse scattering problem is motivated by the
generic setup in imaging, where waves propagate in a reference medium
with smooth wave speed $c(\bx)$, and the goal is to determine
rough perturbations of the medium, the ``reflectivity''. The setup reflects
the separation of scales in the problem, where $c(\bx)$ determines the
kinematics (travel time) of the waves, whereas scattering occurs at
the rough  variations in the medium \cite{symes2009seismic},
like boundaries of inclusions. 

In
some applications, such as  ultrasonic non-destructive evaluation
\cite{drinkwater2006ultrasonic} or  radar imaging
\cite{cheney2009fundamentals}, the reference medium is approximately
homogeneous, like air, so $c(\bx)$ is constant and known. In other
applications, like reflection seismology
\cite{symes2009seismic,Biondi}, $c(\bx)$ must be determined from the
measurements. Velocity estimation is a difficult problem because the
wave fields are oscillatory in time and small perturbations of $c(\bx)$ can
result in travel time perturbations that exceed the period of
oscillations, which is a major change of the wave. This is a serious
issue for data fitting optimization methods that use successive
linearizations, but there are effective approaches for estimating
$c(\bx)$ \cite{stolk2002smooth,uhlmann2001travel,liu1995migration}.

In this paper we assume that $c(\bx)$ is known, and focus attention on
imaging the reflectivity.  Most of the imaging technology is based on
the linearization (Born approximation) of the mapping of the
reflectivity to the scattered wave
\cite{cheney2009fundamentals,beylkin1985imaging,rakesh1988linearised,
  beylkin1990linearized}.  The so-called Kirchhoff formulas
\cite[Chapter 6]{symes1995mathematics} show that if the aperture of
the array is not too large, the Born approximation of the reflected waves  depends to
leading order only on the perturbations of the acoustic impedance
$\sigma(\bx)$. This is the unknown in our setting.

While the linearization assumption has lead to popular imaging methods
known as Kirchhoff migration \cite{Biondi}, matched filtering
\cite{therrien1992discrete} or filtered back-projection
\cite{cheney2009fundamentals}, multiple scattering effects are present
and may lead to significant image artifacts
\cite{delprat2005fundamental,malcolm2007identification}. There has
been progress in the removal of multiple scattering effects in three
different contexts:
\begin{itemize}
\item[(1)] For imaging point-like scatterers buried deep in media with
  small random variations of the wave speed on scales comparable to
  the wavelength
  \cite{borcea2010filtering,borcea2011adaptive,alonso2011detection,
    aubry2009detection}, and for imaging in strongly scattering
  layered media \cite{borcea2012filtering,fomel2007poststack}.
\item[(2)] For imaging scattering surfaces in a smooth reference
  medium, mostly in the context of reflection seismology 
  \cite{weglein1999multiple,malcolm2007identification}  and a related setting   in optics \cite{moskow2012inverse}.
\item[(3)] For imaging almost layered media using the so-called
  Marchenko redatuming method \cite{wapenaar2014marchenko}, and also
  for imaging based on data-driven reduced order models
  \cite{druskin2016direct,DruskinMamonovZaslavsky2017}. This latter
  work is the foundation of the algorithm in this paper.
\end{itemize}

Here we consider an arbitrary unknown acoustic impedance $\sigma(\bx)$
and seek to transform the reflection data \eqref{eq:in5} to
measurements expected in the Born approximation.  The transformation,
called Data to Born (DtB) mapping, is the main result of the paper. We
define it using a reduced order model (ROM) of the wave problem, which
can be calculated from the measurements \eqref{eq:in5} for $q(\bx) =
\ln \sigma(\bx)$.
The ROM is defined by a matrix $\tbL_q$ of special structure,
constructed from the matching relations
\begin{equation}
  \bD_j = \tilde \bM_j(\tbL_q \tbL_q^T), \quad j = 0, \ldots, 2n-1,
  \label{eq:in6}
\end{equation}
for continuously differentiable ROM measurement functions $\tilde
\bM_j$ that do not depend on $q$. These are consistent
with the functions $\bM_j$ in \eqref{eq:in5}, as explained in
\cite{druskin2016direct} and the next sections.

The ROM construction is rooted in the theory of Stieltjes strings due
to Krein \cite{kac1974spectral}.  An outgrowth of this theory, the
spectrally matched grids, also called optimal grids, designed to give
spectrally accurate finite difference approximations of
Dirichlet-to-Neumann maps \cite{druskin1999gaussian}, were used for
discretizations of exterior and multi-scale problems in
\cite{druskin2016near,DMZmultiscale2017}, and for the numerical
solution of the electrical impedance tomography problem in the model
reduction framework in \cite{borcea2011resistor,borcea2014model}. A
related approach, based on Krein's work and the theory of Marchenko,
Gel'fand and Levitan
\cite{krein1951solution,krein1953transfer,marchenko1950,gelfand}, has
been used in inverse hyperbolic problems in layered media in
\cite{gopinath1971inversion,habashy1991generalized,burridge1980gelfand,
symes1979inverse,santosa1982numerical,bube1983one}. Recent
extensions to higher dimensions can be found in
\cite{kabanikhin2011numerical,wapenaar2013three}. At the core of this
theory is the reduction of the inverse scattering problem to a
nonlinear Volterra integral equation, or a system of equations. In the
discrete, linear algebra setting, this translates to the Lanczos and
block Lanczos algorithms or, alternatively, the Stieltjes moment problems
\cite{gallivan1996some,gallivan1996rational,dyukarev2004indeterminacy}
and the Cholesky or block-Cholesky algorithms used in
the construction of $\tbL_q$
\cite{druskin2016direct,DruskinMamonovZaslavsky2017}.

We explain in sections \ref{sect:1D} and \ref{sect:MD} that the matrix
$\tbL_q$ obtained from \eqref{eq:in6} is a Galerkin-Petrov
approximation of the operator $L_q$, for carefully constructed bases
of the spaces of approximation of the fields $\cp(t,\bx)$ and
$\bcu(t,\bx)$. We also discuss in section \ref{sect:1DSPECT} a related
ROM, constructed from spectral measurements of the operator $L_q
L_q^T$ in one dimension \cite{borcea2005continuum}. The analogy is
useful for interpreting the entries of $\tbL_q$ in terms of averages
of the unknown impedance $\sigma(\bx)$ on a special "spectrally matched" grid.

While there are other choices of reduced order models, the ones
considered in this paper have an important property: They are
approximately linear in the unknown $q(\bx)$. This means that if we
had a perturbation $q^\ep(\bx)$ of a known $q^0(\bx)$, of the form
\begin{equation}
  q^\ep(\bx) = q^0(\bx) + \ep \big[q(\bx)-q^0(\bx)\big], \quad 0 < \ep
  \ll 1,
  \label{eq:in7}
\end{equation}
the operator $L_{q}$, which is linear in $q$, would be perturbed as 
\begin{equation}
L_{q^\ep} = L_{q^0} + \ep \big[L_{q} -
    L_{q^0}\big],
  \label{eq:in8p}
\end{equation}
and the corresponding ROM would satisfy a similar relation
\begin{equation}
  \tbL_{q^\ep} \approx \tbL_{q^0} + \ep \big[\tbL_{q} -
    \tbL_{q^0}\big].
  \label{eq:in8}
\end{equation}
Here $\tbL_{q^0}$ is constructed the same way as $\tbL_q$, from the
reference data $\bD^0 = \big\{ \bD_j^0 \big\}_{j=0}^{2n-1}$ calculated
by solving equations \eqref{eq:in.1} with the operator $L_{q^0}$.

We do not have access to the data $\bD_j^\ep$ for coefficient
\eqref{eq:in7}. However, since the ROM is obtained from the matching
conditions \eqref{eq:in6}, we obtain from \eqref{eq:in8} the
approximation \begin{equation} \bD_j^\ep = \tilde\bM_j\big(
  \tbL_{q^\ep} \tbL_{q^\ep}^T \big) \approx \tilde\bM_j\Big(
  \big[\tbL_{q^0} + \ep (\tbL_{q} - \tbL_{q^0})\big] \big[\tbL_{q^0} +
    \ep (\tbL_{q} - \tbL_{q^0})\big]^T \Big).
      \label{eq:in10}
\end{equation}
 The Born data model is defined by
\begin{equation}
  \bD^{Born,\ep}_j = \bD_j^0 + \ep \left[\frac{d}{d \ep'} \bD_j^{\ep'} \big|_{\ep' = 0}\right],
  \label{eq:in11}
\end{equation}
and using \eqref{eq:in10} we approximate it with the DtB mapping
$\mathscr{D}$, which takes the measurements $\bD =\big\{ \bD_j
\big\}_{j=0}^{2n-1}$ with entries \eqref{eq:in6} for the unknown
$q(\bx)$, and returns
\begin{equation}
  \mathscr{D}[\bD] = \Big\{\bD_j^0 + \ep\left[ \frac{d}{d \ep'} \tilde \bM_j \Big(
  \big[\tbL_{q^0} + \ep' (\tbL_{q} - \tbL_{q^0})\big] \big[\tbL_{q^0} +
    \ep' (\tbL_{q} - \tbL_{q^0})\big]^T \Big) \Big|_{\ep' = 0}\right]\Big\}_{j=0}^{2n-1}.
  \label{eq:in12}
\end{equation}
Note that $\ep$ is an arbitrary scaling factor in this equation. We take it 
equal to $1$ so that \eqref{eq:in7} equals $q$.

The DtB mapping \eqref{eq:in12} is described in sections \ref{sect:1D} and
\ref{sect:MD}. We define it form first principles in the one
dimensional case $d=1$ in section \ref{sect:1D}, and then extend the
results to multi-dimensions in section \ref{sect:MD}. The related
inverse spectral problem for the hyperbolic system \eqref{eq:in.1} is
discussed in section \ref{sect:1DSPECT}.  
We end with a summary in
section \ref{sect:sum}.

\section{The DtB mapping in one dimension}
\label{sect:1D}

We  define here  the mapping \eqref{eq:in12} in one
dimension.  We begin in section \ref{sect:1D.1} with the derivation of
the data model \eqref{eq:in5}, starting from the acoustic wave
equation. Then we introduce in section \ref{sect:1D.1p} the wave
propagator operator, which we use in section \ref{sect:1D.2} to
construct the ROM. The matrix $\tbL_q$ that defines the ROM is a
Galekin-Petrov approximation of the operator $L_q$, as shown in
section \ref{sect:1D.3}. The algorithm for computing the DtB map
\eqref{eq:in12} is in section \ref{sect:1D.4}, and we illustrate its
performance with numerical simulations in section \ref{sect:1D.5}.

\subsection{Derivation of the data model}
\label{sect:1D.1}

Let us consider sound waves modeled by the excess acoustic pressure denoted by 
$\mathfrak{p}(t,x)$. We use the different script notation $\mathfrak{p}$ to distinguish this field 
from  another, related pressure field defined below, in equation \eqref{eq:1D14}. 

The pressure 
field $\mathfrak{p}(t,x)$ is defined 
in the domain $x>0$, with sound hard boundary at
$x = 0$, 
\begin{equation}\partial_x \mathfrak{p}(t,0) = 0.
\label{eq:BD1}
\end{equation} 
For a finite
duration $t < t_{max}$, with
\begin{equation}
t_{max} < T_\ell = \int_0^\ell\frac{dx}{c(x)}, 
\label{eq:1D2}
\end{equation}
we can truncate the domain at $x = \ell$ without affecting the wave
at $x = 0$, and set 
\begin{equation}
\mathfrak{p}(t,\ell) = 0.
\label{eq:BD2}\end{equation}
Thus,  $\mathfrak{p}(t,x)$ satisfies the wave equation
\begin{equation}
\big( \partial_t^2 + A \big) \mathfrak{p}(t,x) = \partial_t f(t) \delta(x-0^+),
\quad t \in \mathbb{R}, ~ ~ x \in (0,\ell),
\label{eq:1D4}
\end{equation}
in the domain $\Omega = (0,\ell)$, with boundary condition \eqref{eq:BD1} at the accessible boundary
$\partial \Omega_a = \{0\}$ and \eqref{eq:BD2} at the inaccessible boundary $\partial
\Omega_i = \{\ell\}$. The operator  $A$ is given by 
\begin{equation}
A = - \sigma(\bx) c(x) \partial_x \Big[ \frac{c(x)}{\sigma(x)}
  \partial_x \Big].
\label{eq:1D5}
\end{equation}
The medium is at equilibrium 
prior to the emission of the pulse $f(t)$ from a source located at
$x = 0^+$, 
\begin{equation}
\mathfrak{p}(t,x) = 0, \quad t \ll 0.
\end{equation}

For convenience in the derivation of the ROM, we take
$f(t)$ real valued, with Fourier transform $\hat f(\om) \ge 0$.  For
example, $f(t)$ may be a modulated Gaussian with central frequency
$\om_o$ and bandwidth $B$
\begin{equation*}
f(t) = \frac{\cos (\om_o t)}{\sqrt{2 \pi} B} e^{- \frac{(B t)^2}{2}},
\end{equation*}
so that its Fourier transform is
\begin{equation*}
\hat f(\om) = \int_{-\infty}^\infty dt \, e^{i \om t} f(t) = \frac{1}{2} \Big[ e^{-\frac{(\om-\om_o)^2}{2 B^2}}
  + e^{-\frac{(\om+\om_o)^2}{2 B^2}}\Big].
\end{equation*}

Note that $A$ is self-adjoint in the Hilbert space
$\mathcal{H}_{\frac{1}{\sigma c}} = L^2\Big(
[0,\ell],\frac{1}{\sigma(x) c(x)} dx \Big)$ with weighted inner
product
\begin{equation}
\left<\varphi,\psi\right>_{\frac{1}{\sigma c}} = \int_0^\ell dx \,
\frac{\varphi(x) \psi(x)}{\sigma(x) c(x)}, \quad \forall \varphi, \psi
\in \mathcal{H}_{\frac{1}{\sigma c}},
\label{eq:1D6}
\end{equation}
on the domain of functions $ \varphi(x) \in \mathcal{H}_{\frac{1}{\sigma c}}$, satisfying 
$ \varphi'(0) = 0, ~ \varphi(\ell) = 0 .  $ It has simple and
positive eigenvalues $\{\lambda_j \}_{j \ge 1}$ and the eigenfunctions
$\{y_j(x)\}_{j\ge 1}$ form an orthonormal basis of
$\mathcal{H}_{\frac{1}{\sigma c}}$.  Expanding $\mathfrak{p}(t,x)$ in this basis
we obtain the separation of variables formula
\begin{equation}
\mathfrak{p}(t,x) = f(t) \star \sum_{j=1}^\infty H(t) \cos (t \sqrt{\la_j})
\frac{y_j(0)y_j(x)}{\sigma(0)c(0)},
\label{eq:1D7}
\end{equation} 
where $\star$ denotes convolution, ${H}(t)$ is the Heaviside step
function, and the series is the causal Green's function of 
\eqref{eq:1D4}. 

We work with the even time extension of
$\mathfrak{p}(t,x)$,
\begin{equation}
\mathfrak{p}^{e}(t,x) = \mathfrak{p}(t,x) + \mathfrak{p}(-t,x),
\label{eq:1D9}
\end{equation}
because it has a simpler expression than \eqref{eq:1D7},
\begin{equation}
\mathfrak{p}^e(t,x) = \cos \big(t \sqrt{A}\big) \hat f(\sqrt{A}) \delta(x) =
\sum_{j=1}^\infty\cos (t \sqrt{\la_j}) \hat f(\sqrt{\lambda_j})
\frac{y_j(0)y_j(x)}{\sigma(0)c(0)}.
\label{eq:1D10}
\end{equation}
This defines the data
\begin{equation}
D_j = \mathfrak{p}^e(t_j,0) = \mathfrak{p}(t_j,0) +
\mathfrak{p}(-t_j,0), \quad j = 0, \ldots, 2 n -1,
\label{eq:1D11}
\end{equation}
for the inverse scattering problem with unknown impedance
$\sigma(x)$. The instances $t_j = j \tau$ of measurement are equally
spaced, at sufficiently small interval $\tau = t_{max}/(2 n-1)$, as explained
in the next section. Since $\mathfrak{p}(-t,0) = 0$ for $t$ exceeding the
temporal support of the pulse $f(t)$, the second term in
\eqref{eq:1D11} plays a role only for the first few indexes $j$.

Using the expression \eqref{eq:1D10} and the self-adjointness of $A$,
we can rewrite \eqref{eq:1D11} in the symmetric form
\begin{align}
D_j &= \sigma(0) c(0) \left< \delta(\cdot),\mathfrak{p}^e(t,\cdot)
\right>_{\frac{1}{\sigma c}} = \left< \sqrt{\sigma} \cb,\cos(t_j
\sqrt{A}) \sqrt{\sigma} \cb \right>_{\frac{1}{\sigma
    c}}, \label{eq:1D12}
\end{align}
with the notation 
\begin{equation}
\sqrt{\sigma(x)} \cb(x) = \sqrt{\sigma(0) c(0)} \big[\hat{f}
  \big(\sqrt{A}\big)\big]^{1/2} \delta(x) = \sum_{j=1}^\infty
\big[\hat f(\sqrt{\lambda_j})\big]^{1/2}
\frac{y_j(0)y_j(x)}{\sqrt{\sigma(0)c(0)}}.
\label{eq:1D13}
\end{equation}
We call $b(x)$ the ``sensor function'', because it is supported neat
$x = 0$ and appears in equation \eqref{eq:1D12} as a model of the
source and receiver\footnote{Our construction of the DtB map uses that
  $b(x)$ is supported near $x = 0$, but does not require knowing
  $b(x)$.  In the case of a homogeneous medium we can calculate $b(x)$ in terms of $f(x/c)$, which
  is localized at $x = 0$.  In a variable medium the eigenfunctions
  are not known, but they are oscillatory, and the right hand side in
  \eqref{eq:1D13} is a generalized Fourier series of the smooth
  function $(\hat f \,)^{1/2}$. This series is localized near $x =
  0$.}.

To arrive at the first order hyperbolic system \eqref{eq:in.1}, note that 
\begin{equation}
p(t,x) = \cos(t \sqrt{A}) \sqrt{\sigma(x)} \cb(x)
\label{eq:1D14}
\end{equation}
is the pressure field in the acoustic system of equations
\begin{align}
\partial_t \begin{pmatrix} p(t,x) \\ -u(t,x) \end{pmatrix} = 
\begin{pmatrix} 0 & \sigma(x) c(x) \partial_x \\
\frac{c(x)}{\sigma(x)} \partial_x & 0 \end{pmatrix} \begin{pmatrix}
  p(t,x) \\-u(t,x) \end{pmatrix}, \quad t > 0, ~ ~ x \in (0,\ell),
\label{eq:1D15}
\end{align}
with initial conditions 
\begin{equation} p(0,x) = \sqrt{\sigma(x)} \cb(x), \qquad 
{u}(0,x) = 0,
\end{equation}
and with boundary conditions \begin{equation} {p}(t,\ell) = 0, \qquad 
{u}(t,0) = 0.
\end{equation} 
Here ${u}(t,x)$ is the particle velocity.

The system \eqref{eq:1D15} is not in the desired form for our purpose, because the unknown impedance
$\sigma(x)$ appears in a nonlinear fashion in the coefficients of the differential operator. 
We show next how to transform \eqref{eq:1D15} to the system 
\eqref{eq:in.1}, with operator $L_q$ and its adjoint
$L_q^T$ depending linearly on $q(x) = \ln \sigma(x)$.

\subsubsection{The Schr\"{o}dinger system of equations}

Consider the Liouville transformation
\begin{equation}
\cp(t,x) = \frac{p(t,x)}{\sqrt{\sigma(x)}}, \quad \cu(t,x) =
-\sqrt{\sigma(x)}{u}(t,x),
\label{eq:1D17}
\end{equation}
which takes \eqref{eq:1D15} to 
\begin{align}
\partial_t \begin{pmatrix} \cp(t,x) \\ \cu(t,x) \end{pmatrix} =
\begin{pmatrix} 0 & - L_q \\ L_q^T & 0 \end{pmatrix}
\begin{pmatrix} \cp(t,x) \\ \cu(t,\bx) \end{pmatrix}, \quad t > 0, ~ ~ x
\in (0,\ell).
\label{eq:in.1N}
\end{align}
This is the system \eqref{eq:in.1} in the introduction, with 
\begin{equation}
L_q = - c(x) \partial_x + \frac{1}{2} c(x) \partial_x q(x).
\label{eq:1D18}
\end{equation}
The adjoint of \eqref{eq:1D18} with respect to the inner product
$\left< \cdot, \cdot \right>_{\frac{1}{c}}$ weighted by $1/c(x)$, is
given by
\begin{equation}
L_q^T =  c(x) \partial_x + \frac{1}{2} c(x) \partial_x q(x), 
\label{eq:1D20}
\end{equation}
and we note that both $L_q$ and $L_q^T$ are first order Schr\"{o}dinger operators with potentials that 
are linear in $ q(x) = \ln
\sigma(x).  $

The transformed fields \eqref{eq:1D17}, called henceforth the
``primary wave'' $\cp(t,x)$ and the ``dual wave'' $\cu(t,x)$, satisfy
the initial conditions 
\begin{equation}
 \cp(0,x) = \cb(x), ~ ~  \cu(0,x) = 0,
\label{eq:in.2N}
\end{equation}
and boundary conditions
\begin{equation}
\cu(t,0) = 0, \quad \cp(t,\ell) = 0.
\label{eq:in.3N}
\end{equation}
These are the conditions \eqref{eq:in.3}-\eqref{eq:in.2} stated in the
introduction.

The data
model follows from equations \eqref{eq:1D12}, \eqref{eq:1D14} and
\eqref{eq:1D17}
\begin{equation}
D_j = \left< \cb,\cp(t_j,\cdot)\right>_{\frac{1}{c(x)}} = \int_0^\ell
\frac{dx}{c(x)} \cb(x) \cp(t_j,x), \quad j = 0, \ldots, 2n-1.
\label{eq:1D21}
\end{equation}

\subsubsection{Travel time coordinates}

In one dimension we can avoid dealing with weighted inner products, by
changing coordinates in \eqref{eq:in.1N} from $x$ to the travel time
\begin{equation}
T(x) = \int_0^x \frac{ds}{c(s)}, \quad x \in [0,\ell].
\label{eq:1D23}
\end{equation}
This transformation is invertible for $T \in [0,T_\ell]$, with
$T_\ell$ defined in \eqref{eq:1D2} as the travel time from the
accessible boundary at $x = 0$ to the inaccessible boundary at $x =
\ell$.  Thus, we can write $x = x(T)$, for $T \in [0,T_\ell]$.

We keep the same notation for the operator \eqref{eq:1D18} in the
travel time coordinates
\begin{equation}
L_q = - \partial_T + \frac{1}{2} \partial_T \it{q}(T),
\label{eq:1D24}
\end{equation}
and its adjoint with respect to the usual, Euclidean inner product
$\left<\cdot, \cdot \right>$
\begin{equation}
L_q^T = \partial_T + \frac{1}{2} \partial_T q(T).
\label{eq:1D25}
\end{equation}
We also let $ q(T) = \ln \sigma\big(x(T)\big)  $ and $
\cb(T)=b\big(x(T)\big).  $

The measurements \eqref{eq:1D21} are defined by the primary wave at
the time instances $t_j = j\tau$, denoted by
\begin{equation}
\cp_j(T)= \cp(t_j,x(T)\big) = \cos\Big(t_j \sqrt{L_q L_q^T}\Big)b(T),
\quad j = 0, \ldots, 2n-1.
\label{eq:1D22}
\end{equation}
We also use in the discussion the dual field $\cu(t,x(T))$ evaluated
at the shifted time instances $t_j + \tau/2$. Solving
equations \eqref{eq:in.1N} for $\cu$ we get
\begin{equation} 
\cu\big(t,x(T)) = \sin \Big(t \sqrt{L_q^T L_q}\Big) (L_q^T L_q)^{-1/2}
L_q^T \cb(T),
\label{eq:1D34}
\end{equation}
and we denote
\begin{equation}
\cu_j(T)= \cu(t_j + \tau/2,x(T)\big), \quad j =0,\ldots, 2n-1.
\label{eq:1D34p}
\end{equation}

\subsection{The propagator}
\label{sect:1D.1p}

The propagator of the primary wave is the operator
\begin{equation}
\msP = \cos \Big( \tau \sqrt{L_q L_q^T} \Big),
\label{eq:1D29}
\end{equation}
that maps the initial condition $b(T)$ to $\cp(\tau,x(T))$. We use it
in equation \eqref{eq:1D22} to write 
\begin{equation}
\cp_j(T) = \cos \big( j \arccos (\msP) \big) \cb(T) =
\mathcal{T}_j(\msP) \cb(T),  \quad j = 0, \ldots, 2n-1,
\label{eq:1D30}
\end{equation}
where $ \mathcal{T}_j(\msP) = \cos \big( j \arccos (\msP) \big) $ are
the Chebyshev polynomials of the first kind \cite{rivlin}. The data
model \eqref{eq:1D21} takes the form \eqref{eq:in5}, with measurement
functions $M_j$ defined by\footnote{Note that in our formulation the
  sensor function $b(T)$ depends on $q$. We do not write this
  dependence explicitly in $M_j$ because in
  the ROM construction given in section \ref{sect:1D.2}, $b(T)$ is
  mapped to the "ROM sensor vector" $\tilde \bb = D_0^{1/2} {\bf
    e}_1$, with ${\bf e}_1 = (1, 0,\ldots, 0)^T \in \mathbb{R}^n$ and
  $n \ll N$.  Thus, we can remove the dependence on $q$ of the ROM
  measurement functions $\tilde M_j$ by either normalizing the
  measurements with $D_0$, or by assuming that $q(T)$ is known near
  the accessible boundary i.e., at $T \approx 0$.}
\begin{equation}
D_j = M_j(L_q L_q^T) =\left< \cb,\cp_j\right> = \left< \cb,
\mathcal{T}_j(\msP) \cb\right>, \quad j  = 0, \ldots, 2n-1.
\label{eq:1D33}
\end{equation}

The propagator of the dual wave is the operator
\begin{equation}
\mathscr{U} = \cos \Big(\tau \sqrt{L_q^T L_q}\Big),
\label{eq:1D37}
\end{equation}
and it is shown in \cite[Lemma 3.6]{druskin2016direct} that
\begin{equation}
\cu_j(T) = \Big[ \mathcal{T}_{j}^{(2)}(\mathscr{U}) +
  \mathcal{T}_{j-1}^{(2)}(\mathscr{U}) \Big] \cu_0(T), \quad j = 0, \ldots, 2n-1,
\label{eq:1D36}
\end{equation} 
with $\cu_0$ obtained from \eqref{eq:1D34} evaluated at $t = \tau/2$,
and $\mathcal{T}_j^{(2)}$ the Chebyshev polynomials of the
second kind \cite{rivlin}. 

\subsubsection{Time stepping and factorization of the propagator}

Because the Chebyshev polynomials satisfy the three term recurrence
relation
\begin{align}
&\mathcal{T}_{j+1}(\msP) + \mathcal{T}_{j-1}(\msP) = 2 \msP \,
  \mathcal{T}_j(\msP), \quad j \ge 1, ~ \mbox{and}~ ~ \mathcal{T}_0
  (\msP) = I, \quad \mathcal{T}_1 (\msP) = \msP,
\label{eq:1D31}
\end{align}
where $I$ is the identity operator, we obtain from definition
\eqref{eq:1D30} that the primary wave satisfies the exact time
stepping scheme
\begin{equation}
\frac{1}{\tau^2} \Big[\cp_{j+1}(T)-2\cp_j(T) + \cp_{j-1}(T)\Big] =
-\xi(\msP) \cp_j(T), \quad j = 0,\ldots, 2n-1,
\label{eq:1D39}
\end{equation}
with initial conditions
\begin{equation}
\cp_0(T) = \cb(T), \quad \cp_{-1}(T) = \cp_1(T) = \msP \cb(T).
\label{eq:1D40}
\end{equation}
Here we introduced the affine function
\begin{equation}
\xi(\msP) = \frac{2}{\tau^2}
\big(I - \msP \big),
\label{eq:defXi}
\end{equation}
and the last relation in \eqref{eq:1D40} is derived from
\[
\cp_1(T) + \cp_{-1}(T) = \Big[ \mathcal{T}_{1}(\msP) +
  \mathcal{T}_{-1}(\msP) \Big] \cb(T) = 2 \msP
\cb(T) = 2 \cp_1(T).
\]

Similarly, we obtain an explicit time stepping scheme for the dual
wave, from equation \eqref{eq:1D36} and the definition of the
Chebyshev polynomials of the second kind
\begin{equation}
\mathcal{T}_j^{(2)}(\mathscr{U}) = \left\{ 
\begin{array}{ll} 
\displaystyle 2 \sum_{i=1, i = {\rm odd}}^j
\mathcal{T}_i(\mathscr{U}), \quad &j = {\rm odd}, \\ \displaystyle 2
\sum_{i=0, i = {\rm even}}^j \mathcal{T}_i(\mathscr{U})-I, \quad &j =
    {\rm even},
\end{array} \right.
\label{eq:TCEB1}
\end{equation}
with $\mathcal{T}_{-2}^{(2)}(\mathscr{U}) =
\mathcal{T}_{-1}^{(2)}(\mathscr{U}) = 0.$ We have
\begin{equation}
\frac{1}{\tau^2} \Big[\cu_{j+1}(T)-2\cu_j(T) + \cu_{j-1}(T)\Big] = -\xi(
\mathscr{U} ) \cu_j(T), \quad j = 0,\ldots, 2n-1,
\label{eq:1D39d}
\end{equation}
with $\cu_0(T)$ defined by \eqref{eq:1D34} at $t = \tau/2$ and
\begin{equation}
\cu_0(T) +\cu_{-1}(T) = 0,
\label{eq:1D40d}
\end{equation}
derived from \eqref{eq:TCEB1} and definition \eqref{eq:1D36}.

We can write these two schemes in first order system form, by
factorizing the affine function of the propagators in the right hand
side of \eqref{eq:1D39} and \eqref{eq:1D39d}. We obtain that
\begin{equation}
\xi (\msP) = \frac{4}{\tau^2} \sin^2 \Big(\frac{\tau}{2} \sqrt{L_q
  L_q^T}\Big) = \mathcal{L}_q \mathcal{L}_q^T,
\label{eq:1D41}
\end{equation}
and 
\begin{equation}
\xi (\mathscr{U}) = \frac{4}{\tau^2} \sin^2 \Big(\frac{\tau}{2}
\sqrt{L_q^T L_q}\Big)= \mathcal{L}_q^T \mathcal{L}_q,
\label{eq:1D43}
\end{equation}
with operator
\begin{equation}
\mathcal{L}_q = \frac{2}{\tau} L_q (L_q^T L_q)^{-1/2} \sin
\Big(\frac{\tau}{2} \sqrt{L_q^T L_q}\Big), 
\label{eq:1D42}
\end{equation}
and its adjoint $\mathcal{L}_q^T$ with respect to the Euclidean inner
product. Then, equations \eqref{eq:1D39} and \eqref{eq:1D39d} are
equivalent to the first order time stepping scheme
\begin{align}
\frac{\cp_{j+1}(T)-\cp_j(T)}{\tau} &= - \mathcal{L}_q \cu_j(T), \quad j = 0,
\ldots, 2n-2, \nonumber \\
\frac{\cu_{j}(T)-\cu_{j-1}(T)}{\tau} &=  \mathcal{L}_q^T \cp_j(T), \quad j = 1,
\ldots, 2n-1,
\label{eq:1D44}
\end{align}
with initial conditions 
\begin{equation}
\cp_0(T) = \cb(T), \quad \cu_0(T) +\cu_{-1}(T) = 0.
\label{eq:1D45}
\end{equation}
This is the exact time discretization of the  system
\eqref{eq:in.1N}, for time $t \in [0,t_{max}]$ sampled at intervals $\tau$.

\subsection{The reduced order model}
\label{sect:1D.2}

To avoid technical arguments, we work with the discretization of
\eqref{eq:1D44} on a very fine grid in the interval $(0,T_\ell)$, with
$N \gg 1$ equidistant points at spacing $\Delta T = T_{\ell}/N$. Using
a two point finite difference scheme on this grid, we obtain an $N
\times N$ lower bidiagonal matrix $\bL_q$, the discretization of the
Schr\"{o}dinger operator \eqref{eq:1D24}.  The operator
\eqref{eq:1D42} is discretized by
\begin{equation}
\cbL_q = \frac{2}{\tau} \bL_q (\bL_q^T \bL_q)^{-1/2} \sin
\Big(\frac{\tau}{2} \sqrt{\bL_q^T \bL_q}\Big) = \bL_q \left({\itbf
  I}_{_N} - \frac{\tau^2}{24} \bL_q^T \bL_q + \ldots \right),
\label{eq:1D46}
\end{equation}
where $\itbf{I}_{_N}$ is the $N \times N$ identity matrix. Assuming a
small time sampling interval $\tau$, so that
\begin{equation}
\tau < \|\bL_q^T \bL_q\|^{-1/2} = O \left(\Delta T\right),
\label{eq:1D47}
\end{equation}
we obtain
\begin{equation}
\cbL_q \approx \bL_q.
\label{eq:approxL}
\end{equation}
Then, the primary and dual propagator  $N\times N$ matrices 
\begin{equation}
 \bbP = {\bf I}_{_N} - \frac{\tau^2}{2} \cbL_q \cbL_q^T \approx 
{\bf I}_{_N} - \frac{\tau^2}{2} \bL_q \bL_q^T,
\label{eq:approxP}
\end{equation}
and 
\begin{equation}
 \bbU = {\bf I}_{_N} - \frac{\tau^2}{2} \cbL_q^T \cbL_q \approx 
{\bf I}_{_N} - \frac{\tau^2}{2} \bL_q^T \bL_q,
\label{eq:approxPU}
\end{equation}
are approximately tridiagonal. Here we used definition 
\eqref{eq:defXi}, the factorizations \eqref{eq:1D41}, 
\eqref{eq:1D43} and the approximation \eqref{eq:approxL}.

We call the vectors $\bp_j$ and $\bu_j$ in $\mathbb{R}^N$, with
entries approximating $\cp_j(T)$ and $\cu_j(T)$ on the fine grid, the
primary and dual ``solution snapshots''. They evolve from the initial
values $\bp_0 = \bb$ and $\bu_0$ according to the equations
\begin{equation}
\bp_j = \mathcal{T}_j(\bbP) \bb, \quad \bu_j =
\Big[\mathcal{T}_j^{(2)}(\bbU) + \mathcal{T}_{j-1}^{(2)}(\bbU) \Big]
\bu_0,
\label{eq:solSnaps}
\end{equation}
for $j = 0, \ldots, 2n-1$.  The data model \eqref{eq:1D33} becomes
\begin{equation}
D_j = M_j(L_q L_q^T) \approx M_j (\bL_q \bL_q^T) = \itbf{b}^T\bp_j =
\bb^T \mathcal{T}_j ( \bbP ) \itbf{b},
\label{eq:1D48}
\end{equation}
with small error of the approximation, of order $1/N$, for $N \gg
1$. Here $\bb$ is the ``sensor vector'' in $\mathbb{R}^N$ with entries
defined by the values of the sensor function $b(T)$ on the grid,
multiplied by $\sqrt{\Delta T}$, so that
\begin{equation}
D_0 = \bb^T \bb = \left< b, b \right> + O(1/N).\label{eq:bb}
\end{equation}
We neglect henceforth the $O(1/N)$ error and treat \eqref{eq:1D48} and
\eqref{eq:bb} as equalities.

The ROM is defined by the $n \times n$ symmetric and tridiagonal
(Jacobi) matrix $\tilde{\bbP}$, satisfying the data matching
conditions
\begin{equation}
D_j = \tilde \bb^T \ \mathcal{T}_j (\tilde{\bbP} ) \tilde \bb,
\quad \tilde \bb =D_0^{1/2} {\bf e}_1, \quad j = 0, \ldots, 2n-1,
\label{eq:1D50}
\end{equation}
where ${\bf e}_1 = (1, 0,\ldots, 0)^T \in \mathbb{R}^n$ and $n \ll N$.
Comparing \eqref{eq:1D50} with \eqref{eq:1D48}, we note that $\bb$,
which is supported in the first rows, is replaced by the ``ROM sensor
vector'' $\tilde \bb = \|\bb\| {\bf e}_1$, with $\|\bb\| = D_0^{1/2}$
by \eqref{eq:bb}. We refer to \cite{druskin2016direct} for many
details on the propagator $\tilde{\bbP}$. Here it suffices to obtain
its factorization
\begin{equation}
\tilde{\bbP} = {\itbf I}_{_n} - \frac{\tau^2}{2} \tbL_q \tbL_q^T,
\label{eq:1D51}
\end{equation}
with ${\itbf I}_n$ the $n \times n$ identity, and lower bidiagonal
$\tbL_q$. This is the matrix used in the DtB map \eqref{eq:in12}, and
we explain in the next section how to calculate it.

\subsubsection{Projection ROM}
\label{sect:projROM}

It is shown in \cite[Lemma 4.5]{druskin2016direct} that $\tilde{\bbP}$
can be constructed with an orthogonal projection of $\bbP$ on the span
of the first $n$ primary snapshots $\{\bp_j\}_{j=0}^{n-1}$, the range
of the $N \times n$ matrix
\begin{equation}
\bP = \big( \bp_0, \ldots, \bp_{n-1}\big).
\label{eq:1D53}
\end{equation}
By equation \eqref{eq:solSnaps}, this is the Krylov subspace
\begin{equation}
 \mathcal{K}_n(\bb,\bbP) = \mbox{span}\{ \bb, \bbP \bb,
 \ldots,\bbP^{n-1} \bb \} = \mbox{range}(\bP).
\label{eq:1D52}
\end{equation}
It is intuitive that the projection space is determined only by the
first $n$ snapshots. The backscattered wave measured at $x=0$, for $t
\le t_{2n-1}$, cannot propagate farther than $x(t_{n-1})$ in the
medium, before it reflects and turns back to $x = 0$. This means that
we can image up to depth $x(t_{n-1})$, and all the information is
contained in the subspace \eqref{eq:1D52}.

There are many ways to project on $\mathcal{K}_n(\bb,\bbP)$, depending
on the choice of the basis. We use an orthonormal basis $\{\bv_1,
\ldots, \bv_n\}$ that makes the projection
\begin{equation}
\tilde{\bbP} = \bV^T \bbP \bV,
\label{eq:1D57}
\end{equation}
tridiagonal, where $\bV = (\bv_1, \ldots, \bv_n)$ is the $N \times n$
orthogonal matrix in the QR factorization 
\begin{equation}
\bP = \bV \bR,
\label{eq:1D55}
\end{equation}
with $\bR \in \mathbb{R}^{n \times n}$ invertible and upper triangular
\cite{golub}.  Because of this triangular matrix we obtain from
\eqref{eq:1D55} that the basis satisfies the causality relations
\begin{equation}
\bv_j \in \mbox{span}\{\bp_0, \ldots, \bp_{j-1}\}, \quad j = 1,
\ldots, n.
\label{eq:1D56}
\end{equation}
This is important for at least two reasons: First, it ensures that
$\tilde{\bbP}$ is tridiagonal, as shown in appendix \ref{ap:A}. Second,
it concentrates the support of $\bv_j$ near the wavefront, at $T =
t_j$, and makes the matrix $\bV$ almost independent of the unknown
$q$.

There are two ways of explaining this last property of the
vectors $\bv_j$: One way is to start with $\bv_1$ which equals $\bp_o
= \bb$, up to a normalization factor, and recall that $\bb$ is
supported in the first rows, at travel time $T \approx t_0 = 0$. The
support of the second snapshot $\bp_1$ advances by the travel time $T
= t_1$. Since $\bv_2 \in \mbox{span}\{\bv_1, \bp_1\}$ and $\bv_2$ is
orthogonal to $\bv_1$, the entries in $\bv_2$ must be large at the
wavefront $T = t_1$. Arguing this way, with index $j$ increased one by
one, we see that the support of the orthonormal basis follows the
wavefront of the wave. Depending on how oscillatory the pulse is,
there are some reverberations behind the wavefront, but as shown in
the numerical simulations, the entries in $\bv_j$ are larger around
travel times $T = t_{j-1}$. This property is important in our context,
because the travel times are determined by the known wave speed
$c(x)$, and not the unknown impedance or, equivalently, $q(x)$. This
means that $\bV$ is almost independent of $q(x)$, as illustrated in
section \ref{sect:1D.4}.

The other way of explaining is algebraic: By causality, the matrix
$\bP$ of the primary snapshots is approximately upper triangular. The
approximation is because $\bP \in \mathbb{R}^{N \times n}$ is a tall
rectangular matrix and $\bb$ is not an exact delta-function, but an
approximation. If there where no inhomogeneities in the medium, there
would be no reflected waves and the matrix $\bP$ would be
approximately diagonal. The inhomogeneities cause reflections, which
fill-in the upper triangular part of $\bP$. The QR orthogonalization
\eqref{eq:1D55} transforms the {almost} upper triangular matrix $\bP$
to the {almost} identity matrix $\bV$, which is almost independent of
$q$.

\subsubsection{The calculation of the ROM}
\label{sect:calcROM}

Although the QR factorization \eqref{eq:1D55} is useful for
understanding the ROM, we cannot use it directly to compute
$\tilde{\bbP}$ because we do not know the matrix \eqref{eq:1D53}. We
only know the inner products of its columns with $\bb$, from
\eqref{eq:1D48}.  We now explain how to calculate $\tilde{\bbP}$ from
the matching relations \eqref{eq:1D50}.

Let us begin with the calculation of the upper triangular matrix $\bR$. We
obtain from equations \eqref{eq:solSnaps} and \eqref{eq:1D55} that
\begin{equation}
\big(\bP^T \bP\big)_{j,k} =  \bb^T
\mathcal{T}_{j-1}(\bbP) \mathcal{T}_{k-1}(\bbP) \bb = (\bR^T \bR)_{j,k}, \quad 
j,k = 1, \ldots, n,
\label{eq:MassM}
\end{equation}
where we used the symmetry of $\bbP$. The
Chebyshev polynomials satisfy the relation
\begin{equation}
\label{eq:ChebT4}
\mathcal{T}_j(\bbP) \mathcal{T}_k(\bbP) = \frac{1}{2} \Big[
  \mathcal{T}_{j+k}(\bbP) + \mathcal{T}_{|j-k|}(\bbP)\Big],
\end{equation}
so substituting in \eqref{eq:MassM} and using \eqref{eq:1D48}, we get
\begin{equation}
\big(\bP^T \bP\big)_{j,k} = (\bR^T \bR)_{j,k} = \frac{1}{2}
\Big(D_{j+k-2} + D_{|j-k|}\Big).
\label{eq:MassM1}
\end{equation}
This shows that the $n \times n$ matrix $\bP^T \bP$ can be determined
from the data, and that $\bR$ can be calculated from its Cholesky
factorization \cite{golub}
\begin{equation}
\bP^T \bP = \bR^T \bR.
\label{eq:Chol1}
\end{equation} 

With the matrix $\bR$ calculated from \eqref{eq:Chol1}, we solve for
$\bV$ in \eqref{eq:1D55} to obtain
\begin{equation}
\bV = \bP \bR^{-1},
\label{eq:1D55p}
\end{equation}
and then rewrite \eqref{eq:1D57} as
\begin{equation}
\tilde{\bbP} = \bR^{-T} (\bP^T \bbP \bP) \bR^{-1}.
\label{eq:1D56p}
\end{equation}
The matrix in parentheses has the entries
\begin{equation}
(\bP^T \bbP \bP)_{j,k} = \bb^T \mathcal{T}_{j-1}(\bbP) \bbP
  \mathcal{T}_{k-1}(\bbP) \bb, 
\label{eq:Stiff1}
\end{equation}
by definition \eqref{eq:solSnaps}. Then, relation \eqref{eq:ChebT4},
$\mathcal{T}_1(\bbP) = \bbP$, and definition \eqref{eq:1D48} give
\begin{equation}
(\bP^T \bbP \bP)_{j,k} = \frac{1}{4} \Big(D_{j+k-1} + D_{|k-j+1|} + 
D_{|k-j-1|} + D_{|k+j-3|} \Big),
\label{eq:Stiff2}
\end{equation}
for $j,k = 1, \ldots, n$. This shows that $\bP^T \bbP \bP$ can be
computed from the data, and the propagator $\tilde{\bbP}$ follows from
\eqref{eq:1D56p}.

To obtain the factorization \eqref{eq:1D51}, we note from
\eqref{eq:1D57} that the spectral norm of the ROM propagator is
bounded above by the spectral norm of $\bbP = \cos\big(\tau
\sqrt{\bL_q^T \bL_q}\big)$.  With our choice \eqref{eq:1D47} of $\tau$
this norm is strictly less than one, so ${\bf I}_{_n} - \tilde{\bbP}$
is positive definite.  Therefore, we can calculate the matrix $\tbL_q$
in \eqref{eq:1D51} from another Cholesky factorization
\begin{equation}
\xi\big(\tilde{\bbP}\big) = \frac{2}{\tau^2} \big({\itbf I}_{_n} -
\tilde{\bbP}\big) = \tbL_q \tbL_q^T.
\label{eq:Chol}
\end{equation}
This is the ROM version of equation \eqref{eq:1D41}.

It remains to show that the vector $\tilde \bb$ in the data matching
conditions has the form given in \eqref{eq:1D50}. We define
$\tilde{\bb}$ as the projection of $\bb$ on the space \eqref{eq:1D52},
given by
\begin{equation}
\tilde \bb = \bV^T \bb. \label{eq:bb1}
\end{equation}
Using equations \eqref{eq:Chol1}, \eqref{eq:1D55p}, \eqref{eq:bb} and the upper 
triangular structure of $\bR$ we get 
\begin{equation}
\tilde \bb = \bR^{-T} \bP^T \bP {\bf e}_1 = \bR {\bf e}_1  = 
\big({\bf e}_1^T \bP^T \bP {\bf e}_1\big)^{1/2} {\bf e}_1 = 
\big(\bb^T \bb\big)^{1/2}{\bf e}_1 = D_0^{1/2}{\bf e}_1,
\label{eq:bb2}
\end{equation}
as stated in \eqref{eq:1D50}.

The ROM measurement functions $\tilde M_j$ in \eqref{eq:in6} are defined by 
\begin{equation}
\tilde M_j(\tbL_q \tbL_q^T) = D_0^{1/2} {\bf e}_1^T \mathcal{T}_j
\Big( {\itbf I}_{_n} - \frac{\tau^2}{2} \tbL_q \tbL_q^T\Big) {\bf e}_1
D_0^{1/2}, \quad j = 0, \ldots, 2n-1.
\label{eq:ROMMeas}
\end{equation}
We can make them independent of $q$ by normalizing the measurements
with $D_0$, which is strictly positive by
\eqref{eq:1D17}-\eqref{eq:in.1N}.  Alternatively, we may suppose that
$q(T)$ is known near $T = 0$, and conclude from the causality of the
wave equation that $D_0$ is independent of the variations of $q(T)$ at
larger $T$.  We make this assumption henceforth, and treat $D_0$ as
constant.

\subsection{The Galerkin-Petrov approximation}
\label{sect:1D.3}

Here we show that the lower bidiagonal matrix $\tbL_q$ computed above
is a Galerkin-Petrov approximation of the operator $\cbL_q$ in
\eqref{eq:1D46}, which in turn is an approximation of the
Schr\"{o}dinger operator $\bL_q$.

Multiplying \eqref{eq:Chol} on the right with the inverse of
$\tbL_q^T$, denoted by $\tbL_{q}^{-T}$, we have
\begin{equation}
\tbL_q = \xi(\tilde{\bbP}) \tbL_{q}^{-T} = \bV^T
\Big[\frac{2}{\tau^2} \big({\itbf I}_{N} - \bbP
  \big)\Big]\bV\tbL_{q}^{-T} = \bV^T \cbL_q \cbL_q^T \bV \tbL_{q}^{-T},
\label{eq:GP1}
\end{equation}
where we used definitions \eqref{eq:approxP}, \eqref{eq:1D57} and
$\bV^T \bV = {\itbf I}_n$.  We rewrite the result as 
\begin{equation}
\tbL_q = \bV^T \cbL_q \bW,
\label{eq:GP2}
\end{equation}
using the $N \times n$ matrix
\begin{equation}
\bW = \cbL_q^{T} \bV \tbL_{q}^{-T},
\label{eq:GP3}
\end{equation}
which is orthogonal by equation \eqref{eq:Chol},
\begin{equation}
\bW^T \bW = \tbL_q^{-1} \bV^T \cbL_q \cbL_q^T \bV \tbL_{q}^{-T} =
\tbL_q^{-1}\xi(\tilde{\bbP}) \tbL_q^{-T} = {\itbf I}_n.
\end{equation} 
Thus, we conclude that $\tbL_q$ is the Galerkin-Petrov approximation
of the operator $\cbL_q$, with the primary field approximated in the
space \eqref{eq:1D52}, and the dual field approximated in the range of
$\bW$. This is the same as the range of the matrix $\bU = (\bu_0,
\ldots,\bu_{n-1})$ of the dual snapshots, as explained in section
\ref{sect:OSN}.

\begin{rem}
\label{rem:Lin}
It follows from \eqref{eq:approxL} and the linearity of $\bL_q$ with
respect to $q$ that $\cbL_q $ is approximately linear in $q$. The
discussion at the end of section \ref{sect:projROM}, which is for the
columns of matrix $\bV$, but extends verbatim to matrix $\bW$,
explains that $\bV$ and $\bW$ are almost independent of $q$. Thus,
equation \eqref{eq:GP2} yields approximate linearity of the reduced
order matrix $\tbL_q$ in $q$, as needed in the DtB mapping.
\end{rem}

\subsection{The data to Born  mapping}
\label{sect:1D.4}

Let $\tbL_{q^0}$ be the ROM matrix in the reference medium with known
impedance $\sigma^0(x)$ and Schr\"{o}dinger potential $q^0(x) = \ln
\sigma^0(x)$. Let also $\bV^0$ and $\bW^0$ be the projection matrices
in this medium. As explained in section \ref{sect:projROM} and Remark
\ref{rem:Lin}, these matrices change slowly with the potential $q(x)$,
so for the perturbed $q^\ep(x)$ defined in \eqref{eq:in7} we have
\begin{equation}
\bV^\ep \approx \bV^0, \quad \bW^\ep \approx \bW^0. 
\label{eq:approxVW}
\end{equation}
Equation \eqref{eq:GP2} gives 
\begin{equation}
\tbL_{q^\ep} \approx \tbL_{q^0} + \bV^{0^T}
\left(\cbL_{q^\ep}-\cbL_{q^0}\right) \bW^0,
\label{eq:DtB1}
\end{equation}
and due to the approximation \eqref{eq:approxL} and the linearity of
$\bL_q$ in $q$, we have 
\begin{equation}
\cbL_{q^\ep} - \cbL_{q^0} \approx \bL_{q^\ep} - \bL_{q^0} = \ep \big(
\bL_q - \bL_{q^0}\big).
\label{eq:DtB2}
\end{equation}
Substituting \eqref{eq:DtB2} in \eqref{eq:DtB1} we get the approximate 
linearity  relation
\eqref{eq:in8}, which makes the mapping \eqref{eq:in12}
useful.

\vspace{0.05in}
\begin{alg}
\label{eq:algo1} 
The algorithm for computing the DtB mapping is as follows:

\vspace{0.05in} \noindent \textbf{Input:} data ${\itbf D} =
\{D_j\}_{j=0}^{2n-1}.$

\vspace{0.05in} \noindent 1. Map the data to the ROM matrix $\tbL_{q}$
using equations \eqref{eq:MassM1}, \eqref{eq:1D56p}, \eqref{eq:Stiff2}
and the Cholesky factorizations \eqref{eq:Chol1} and \eqref{eq:Chol}.

\vspace{0.05in} \noindent 2. Compute the data ${\itbf D}^0 =
\{D_j^0\}_{j=0}^{2n-1}$ in the reference medium with given $q^0$,
using formula
\begin{equation}
D_j^0 = \bb^{0^T} \cos \Big(j\tau \sqrt{\bL_{q^0} \bL_{q^0}^T}\Big)\bb^0, \quad
j = 0, \ldots, 2n-1,
\end{equation}
where $\bb^0$ is the sensor vector in the reference medium and
$\bb^{0^T}$ is its transpose. Moreover,  $\bL_{q^0}$ is the $N \times N$ lower bidiagonal 
matrix, the discretization on the fine grid 
with $N$ points of the operator \eqref{eq:1D24} with reference potential $q^0$.

\vspace{0.05in} \noindent 3. Map the data ${\itbf D}^0$ to the ROM matrix $\tbL_{q^0}$
using equations \eqref{eq:MassM1}, \eqref{eq:1D56p}, \eqref{eq:Stiff2}
and  Cholesky factorizations \eqref{eq:Chol1}, \eqref{eq:Chol}.

\vspace{0.05in} \noindent 4. Use definition \eqref{eq:ROMMeas} and
 the chain rule given by Algorithm \ref{alg:chain} to calculate
\begin{align}
& \frac{d}{d \ep} D_0^{-1/2} \tilde M_j \Big( \tbL_\ep \tbL_\ep^T \Big)
  D_0^{-1/2} \Big|_{\ep = 0} = {\bf e}_1^T \frac{d}{d \ep}
  \mathcal{T}_j \Big( {\itbf I}_{_n} - \frac{\tau^2}{2}\tbL^\ep
  \tbL^{\ep^T}\Big) \Big|_{\ep = 0} {\bf e}_1 \, \label{eq:chain}
\end{align}
where $ \tbL^\ep = \tbL_{q^0} + \ep\big(\tbL_q-\tbL_{q^0}\big) $ and
$\tbL^{\ep^T}$ is its transpose. Recall that $D_0$ is the first
measurement, at  $t_0 = 0$.

\vspace{0.05in}\noindent \textbf{Output:} the transformed data
$\mathscr{D}[{\itbf D}]$ given by \eqref{eq:in12}, with the derivative
calculated from \eqref{eq:chain} and $\ep = 1$.
\end{alg}

\vspace{0.05in} The derivative in \eqref{eq:chain} must be calculated
carefully, because the usual chain rule does not apply to matrix
functions, unless the matrix commutes with its  derivative.  We
obtain it in the next algorithm using the recursion relation
\eqref{eq:1D31}.
\vspace{0.05in}

\begin{alg} The algorithm for computing  \eqref{eq:chain} is as follows:
\label{alg:chain}\\
\textbf{Input:} $\tbL_q, \tbL_{q^0}$ and $\tau$.

\vspace{0.05in} \noindent Denote $\tilde {\bbP^\ep} = {\itbf I}_n -
\frac{\tau^2}{2} \tbL^\ep \tbL^{\ep^T},~$ $ \tilde \bp_{j}^\ep =
\mathcal{T}_j(\tilde {\bbP^\ep}) {\bf e}_1$ and $~\tilde {\bf z}_j =
\frac{d}{d \ep} \tilde \bp_{j}^\ep \Big|_{\ep = 0},$ $~j = 0,
\ldots, 2n-1.$

\vspace{0.05in} \noindent We have $ \tilde {\bf z}_0 = {\bf 0}$ and
$\tilde {\bf z}_1 = - \frac{\tau^2}{2} \Big( \tbL_{q} \tbL_{q^0}^T +
\tbL_{q^0} \tbL_q^T - 2 \tbL_{q^0} \tbL_{q^0}^T \Big) {\bf e}_1, $ and
for $j = 2, \ldots, 2n-1$,
\[
\tilde{\bf z}_j = 2 \xi (\tilde{\bbP^0}) \tilde {\bf z}_{j-1} -
\tilde{\bf z}_{j-2} -\tau^2\Big( \tbL_{q} \tbL_{q^0}^T + \tbL_{q^0}
\tbL_q^T - 2 \tbL_{q^0} \tbL_{q^0}^T \Big) \tilde \bp_{j-1}^0.
\]

\vspace{0.02in} \noindent \textbf{Output:} The derivative
\eqref{eq:chain} given by $ \frac{d}{d \ep} D_0^{-1/2} M_j \Big(
\tbL_\ep \tbL_\ep^T \Big) D_0^{-1/2} \Big|_{\ep = 0} = {\bf e}_1^T
\tilde {\bf z}_j.$

\end{alg}
\subsection{Numerical results}
\label{sect:1D.5}
\begin{figure}
\centering
\begin{tabular}{cc}
\includegraphics[width=0.44\textwidth]{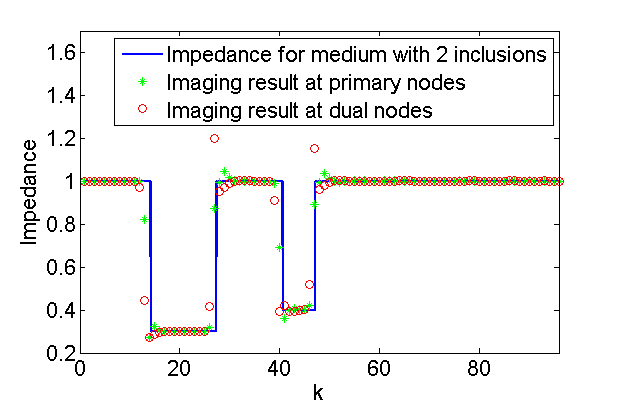} \\
\end{tabular}
\vspace{-0.14in}
\caption{Layered relative acoustic impedance model $\sigma(x)$ and discrete
  inversion results.  The discrete inversion is discussed in
  section~\ref{sect:OSN}. The abscissa is the spatial (primary) grid node index, sampled at
  interval $\tau$. }
\label{fig:1dmodel}
\end{figure}

We present numerical results for a layered model, with relative
acoustic impedance shown in Figure~\ref{fig:1dmodel}.  The relative
impedance is defined as the ratio of the impedance and that of the
homogeneous background. We display it as a function of the travel
time, at steps $\tau$ chosen consistent with the Nyquist sampling rate
of the Gaussian pulse used in the simulations.  To avoid the ``inverse
crime'', the data are generated with a finite-difference time-domain
algorithm, on an equidistant grid with steps much smaller than $\tau$.

\begin{figure}
\centering
\begin{tabular}{cc}
\includegraphics[width=0.5\textwidth]{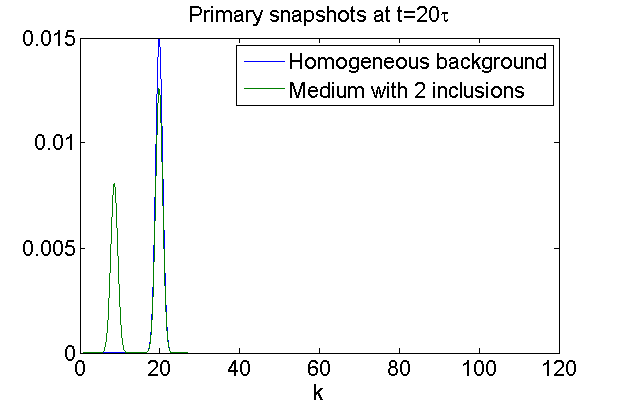} &
\includegraphics[width=0.46\textwidth]{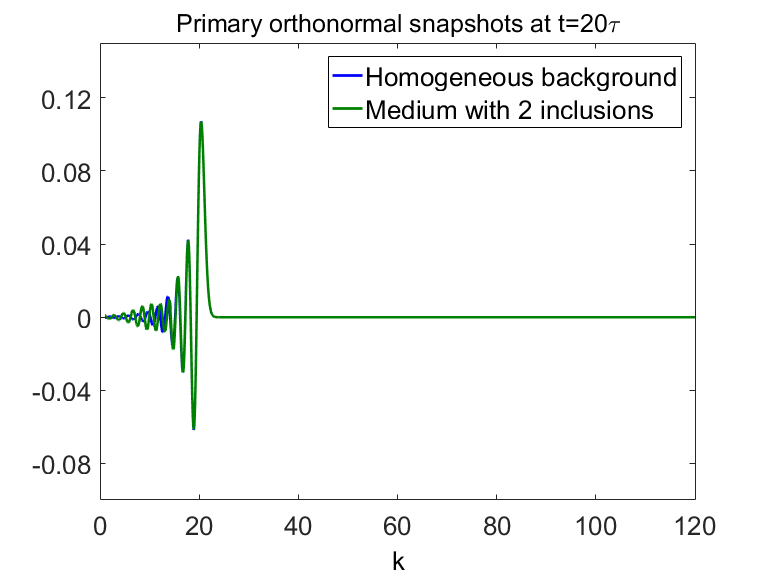}\\
\includegraphics[width=0.5\textwidth]{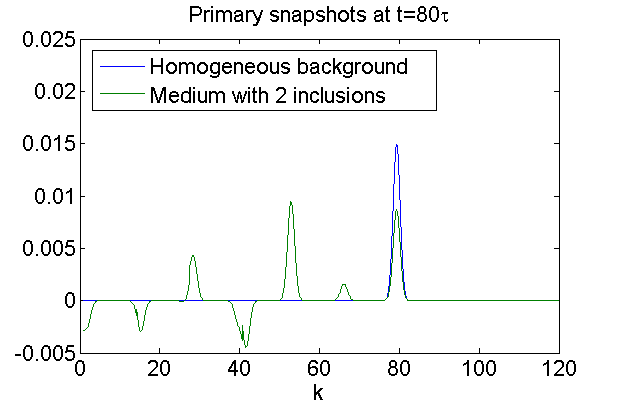} &
\includegraphics[width=0.46\textwidth]{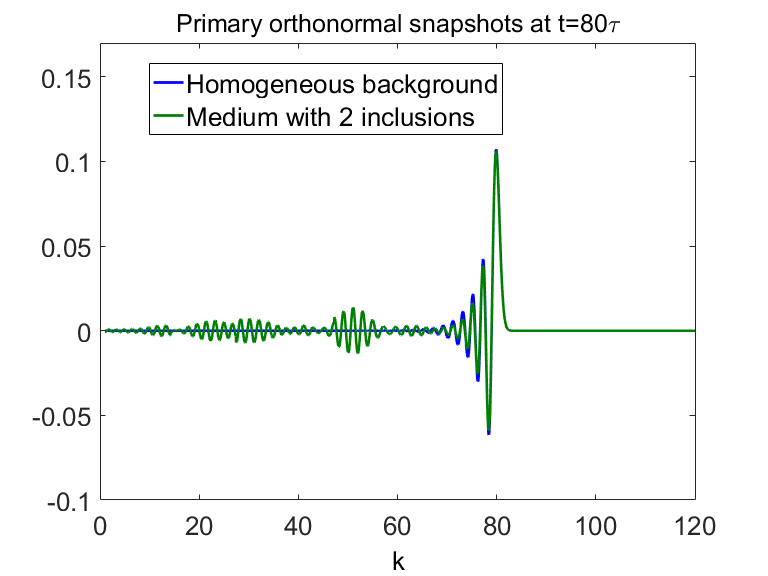} 
\end{tabular}
\vspace{-0.12in}\caption{Primary snapshots $\bp_j$ (left column) and primary
  orthonormal snapshots $\bv_j$ (right column) at time index $j =
  20$ (top row) and $j = 80$ (bottom row). The abscissa is the travel
  time index, sampled at interval $\tau$.}
\label{fig:model&snapshots}
\end{figure}

In Figure~\ref{fig:model&snapshots} we show the primary snapshots,
columns of $\bP^0$ and $\bP$, for the homogeneous background with
relative impedance $\sigma^0 = 1$, and the layered model.  The wave
has crossed all the discontinuities of the impedance by the time $t =
80 \tau$ of the snapshots displayed in the bottom row.  Thus, we
observe significant differences between $\bP_{80}$ and
$\bP^0_{80}$. These consist of the decrease of the amplitude of the
first arrival and the large multiple reflections.

In Figure~\ref{fig:model&snapshots} we show the 
columns of $\bV^0$ and $\bV$, for the homogeneous background and
the layered model.  We call these columns the primary orthonormal snapshots. We observe that they are almost independent
of the medium, as discussed in Remark \ref{rem:Lin}. A similar
behavior holds for the dual orthonormal snapshots, not shown here.

\begin{figure}
\centering
\begin{tabular}{cc}
\includegraphics[width=0.47\textwidth]{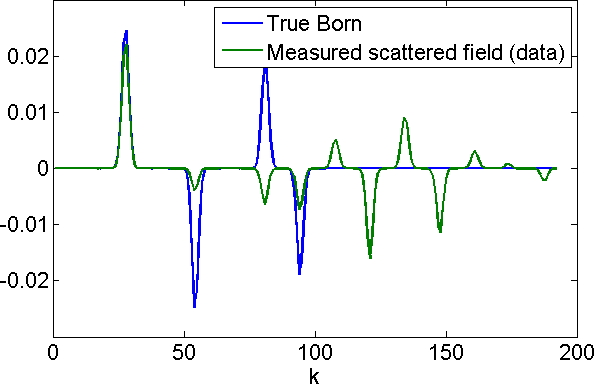} & 
\includegraphics[width=0.45\textwidth]{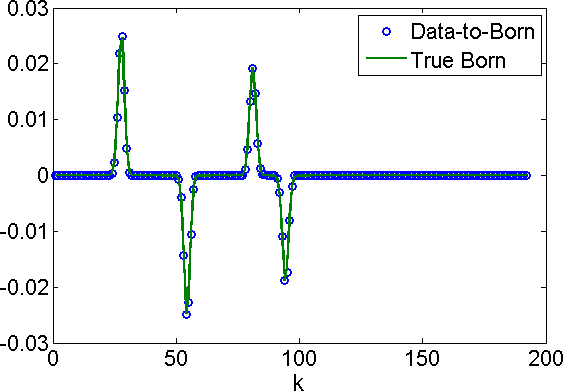} 
\end{tabular}
\vspace{-0.12in}\caption{Left: the raw scattering data (green line) and the Born
  approximation (blue line).  Right: the data transformed with the DtB
  map (blue line) is indistinguishable from the Born approximation
  (green line). The abscissa is the time index.}
\label{fig:dtb1d}
\end{figure}

In Figure~\ref{fig:dtb1d} we show the raw scattering data, its Born
approximation\footnote{The Born approximation cannot be computed in 
the inverse scattering problem, because it corresponds to solving 
the  wave equation linearized with respect to 
the logarithm of the unknown impedance. This is why we need the DtB transform. We display the Born approximation
only 
for comparison with the output of the DtB algorithm.}  and the data obtained with the DtB algorithm. We observe
that the strong multiples in the raw data are removed, and that the
result is indistinguishable from the Born approximation.

\begin{rem}
\label{rem:Num}
Our experiments with different $\tau$ (not showed here) indicate that
the $l_\infty$ discrepancy between the true Born approximation and the
output of the DtB algorithm decays as $O(\tau^2)$ for smooth
$\sigma(x)$, in agreement with the approximation error in
\eqref{eq:1D46}. We speculate that if we solved exactly for $\bL_q$ in
\eqref{eq:1D46}, the discrepancy would decay exponentially in
${\tau}^{-1}$.
\end{rem}
\section{A related inverse spectral problem}
\label{sect:1DSPECT}

In this section we look in more detail at the entries of the ROM
matrix $\tbL_q$, and compare it with another ROM obtained from
spectral measurements of the operator \eqref{eq:1D5} in the wave
equation. 

\subsection{The orthonormal snapshots and the entries in $\tbL_q$} 
\label{sect:OSN}

The first primary snapshots $\{\bp_j\}_{j=0}^{n-1}$ span the Krylov
space $\mathcal{K}_n(\bb,\bbP)$ defined in \eqref{eq:1D52}, and the
dual snapshots $\{\bu_j\}_{j=0}^{n-1}$ span the Krylov space
$\mathcal{K}_n(\bU_0,\bbU)$, as follows from definition
\eqref{eq:solSnaps}. 

The classical method for computing an orthonormal basis of a Krylov
subspace is given by the Lanczos method \cite{parlett1998symmetric},
which is used in \cite[Algorithm 3.1]{druskin2016direct} to calculate
the orthogonal  vectors\footnote{The vectors $\overline{\bf p}_j$ and $\overline{\bf u}_j$  are 
called orthogonalized primary and dual snapshots in 
\cite{druskin2016direct}.}
$\{\overline{\bf p}_j\}_{j=1}^n$ and the orthogonal  vectors
 $\{\overline{\bf u}_j\}_{j=1}^n$,  satisfying the equations
\begin{align}
\frac{\overline{\bf p}_{j+1} - \overline{\bf p}_{j}}{\gamma_j} = - \cbL_q
\overline{\bf u}_j, \quad \frac{\overline{\bf u}_{j} -
  \overline{\bf u}_{j-1}}{\hat \gamma_j} = \cbL_q^T \overline{\bf p}_j,
\label{eq:O1}
\end{align}
for $j \ge 1$, with initial conditions 
$
\overline{\bf p}_{1} = \bb,$ $\overline{\bf u}_{0} = 0,$
and coefficients
\begin{equation}
\gamma_j = \frac{1}{\|\overline{\bf u}_j\|^2}, \quad 
\hat \gamma_j = \frac{1}{\|\overline{\bf p}_j\|^2}, \quad j = 1, \ldots, n.
\label{eq:O3}
\end{equation}
It is also shown in \cite[Sections 4.2, 4.3]{druskin2016direct} how
the coefficients \eqref{eq:O3} enter the expression of the tridiagonal
ROM propagator $\tilde{\bbP}$. Using the factorization
\eqref{eq:Chol}, we obtain from those results the lower bidiagonal
matrix $\tbL_q$ with entries 
\begin{align}
(\tbL_q)_{j,j} &= -1/\sqrt{\gamma_j \hat \gamma_j}, \quad j = 1,
  \ldots, n, \nonumber \\ 
(\tbL_q)_{j+1,j} &= 1/\sqrt{\gamma_{j} \hat
    \gamma_{j+1}}, \quad j = 1, \ldots, n-1.
\label{eq:O3p}
\end{align}

The columns of the projection matrix $\bV$ on the Krylov space
$\mathcal{K}_n(\bb,\bbP)$, aka the orthonormal primary snapshots, are given by 
\begin{equation}
\bv_j = \sqrt{\hat \gamma_j} \, \overline{\bf p}_j, \quad j = 1,
\ldots,n,
\label{eq:O4}
\end{equation}
and the projection matrix $\bW$ satisfies by definition \eqref{eq:GP3}
\begin{equation}
\bW \tbL_q^T = \cbL_q^T \bV.
\label{eq:O5}
\end{equation}
The $j$-th column in this equation reads 
\begin{equation}
\frac{1}{\sqrt{\hat \gamma_j}} \bW \tbL_q^T \itbf{e}_j = \cbL_q^T
\overline{\bf p}_j = \frac{\overline{\bf u}_{j} -
  \overline{\bf u}_{j-1}}{\hat \gamma_j},
\label{eq:O6}
\end{equation}
where we used \eqref{eq:O1} and \eqref{eq:O4}. The left hand side is a linear combination of the columns $\bw_j$ and
$\bw_{j-1}$ of $\bW$, because $\tbL_q$ is lower bidiagonal.  Using the
expression \eqref{eq:O3p} of the entries in $\tbL_q$ and equation
\eqref{eq:O4}, we conclude that
\begin{equation}
\bW = \Big(\bw_1, \ldots, \bw_n), \quad \bw_j = \sqrt{\gamma_j} \,
\overline{\bf u}_j.
\label{eq:O7}
\end{equation}
This shows that $\bW$ is the matrix of orthormal  dual
snapshots, as stated in the previous section.

Let us take the constant reference impedance $\sigma^0= 1$,
corresponding to the potential $ q^0 = 0$, and define the coefficients
\begin{align}
\sigma_j = \frac{\hat \gamma_j^0}{\hat \gamma_j} =
\frac{\|\overline{\bf p}_j\|^2}{\|\overline{\bf p}^0_j\|^2},
\qquad  \hat \sigma_j = \frac{ \gamma_j}{ \gamma_j^0} =
\frac{\|\overline{\bf u}^0_j\|^2}{\|\overline{\bf u}_j\|^2}.
\label{eq:ratiosROM}
\end{align}
With these coefficients we introduce the discrete Liouville transform
\begin{equation}
\overline{\bp}_j =\frac{\overline{\bf p}_j}{\sqrt{\sigma_j}},
\quad \overline{\bu}_j = \sqrt{\hat \sigma_j} \,\overline
{\bf u}_j,
\label{eq:LD1}
\end{equation}
substitute it in \eqref{eq:O1} and obtain after straightforward
algebraic manipulations the system
\begin{align}
\frac{\overline{\bp}_{j+1}- \overline{\bp}_j
}{\gamma_j^0 } + \overline{\bp}_{j+1}
\left(\frac{\sqrt{\sigma_{j+1}}-\sqrt{\hat \sigma_j}}{\gamma_j^0
  \sqrt{\hat \sigma_j}}\right) + \overline{\bp}_j
\left(\frac{\sqrt{\hat \sigma_j}-\sqrt{\sigma_j}}{\gamma_j^0
  \sqrt{\hat \sigma_j}}\right) &= -\pmb{\mathcal{L}}_q
\overline{\bu}_j, \label{eq:LD2}\\ \frac{\overline{\bu}_{j}-
  \overline{\bu}_{j-1} }{\hat \gamma_j^0 } -
\overline{\bu}_{j} \left(\frac{\sqrt{\hat \sigma_{j}}-\sqrt{
      \sigma_j}}{\hat \gamma_j^0 \sqrt{\hat \sigma_j}}\right) -
  \overline{\bu}_{j-1} \left(\frac{\sqrt{
      \sigma_j}-\sqrt{\hat \sigma_{j-1}}}{\hat \gamma_j^0 \sqrt{\hat
      \sigma_{j-1}}}\right) &= \pmb{\mathcal{L}}_q^T
  \overline{\bp}_j, \label{eq:LD3}
\end{align}
for $j = 1, \ldots, n$. Recalling the approximation
\eqref{eq:approxL}, we see that the finite difference operators in the
left hand sides of equations \eqref{eq:LD2}-\eqref{eq:LD3} can be
interpreted as discretizations of $ L_q ^T = \partial_T + \frac{1}{2}
\partial_T q(T) $ and $-L_q = \partial_T -\frac{1}{2} \partial_T q(T)
$, with $q(T) = \ln \sigma(T)$.  The discretization is on a special
grid with primary points spaced at $\gamma_j^0$, and dual points
spaced at $\hat \gamma_j^0$.  In our case these equal $\tau$
\cite{druskin2016direct}.  In Figure~\ref{fig:model&snapshots} the
primary grid corresponds to the integer values in the abscissa, and
the dual grid points are offset by $\tau/2$. We note that the peaks of
the primary orthonormal snapshots $\bV_j$ are approximately aligned with
the $j$-th primary grid point, which is the location of the wavefront.

The coefficients $\sigma_j$ and $\hat \sigma_j$ are approximations of
the impedance at the primary and dual grid points, and the terms in
the parentheses in \eqref{eq:LD2}-\eqref{eq:LD3} are
discretizations of
\[
\frac{1}{2}\partial_T q(T) = \frac{\partial_T
  \sqrt{\sigma(T)}}{\sqrt{\sigma(T)}}.
\]
We display in Figure~\ref{fig:1dmodel} the values $\sigma_j$ and $\hat
\sigma_j$ computed for the impedance model considered in
section~\ref{sect:1D.4}. They give a reasonable approximation of the
discontinuous impedance, but better results can be obtained by
inverting the data processed by the DtB algorithm, which is almost
indistinguishable from the Born approximation.

We show next that the same approximation formulas \eqref{eq:ratiosROM}
arise for another ROM, constructed from different measurement
functions of the operator $L_q L_q^T$ than in \eqref{eq:1D48}. The
values of the ROM coefficients $\{\gamma_j,\hat \gamma_j\}_{j=1}^n$
are different, but the same ratios $\{\hat \gamma_j^0/\hat \gamma_j,
\gamma_j/\gamma_j^0\}_{j=1}^n$ define approximations of $\sigma(T)$ on
the ROM dependent grids with primary and dual point spacings defined
by $\{\gamma_j^0,\hat \gamma_j^0\}_{j=1}^n$. With this ROM described
below, it is proved in \cite{borcea2005continuum} that the
approximations \eqref{eq:ratiosROM} converge to the unknown impedance
function $\sigma(T)$ in the limit $n \to \infty$.

\subsection{The spectrally matched ROM}
\label{sect:SPECT}

In this section we draw an analogy between the ROM constructed from
the data matching conditions \eqref{eq:in6} and the "spectrally
matched" ROM introduced and analyzed in \cite{borcea2005continuum}.
Spectrally matched means that the ROM defines a three point finite
difference scheme in $x$ for the wave equation satisfied by the
pressure field ${p}(t,x)$ in \eqref{eq:1D15}, modeled with an
$n \times n$ tridiagonal matrix constructed from the truncated
spectral measure of the differential operator $A$ in \eqref{eq:1D5}.

The Laplace transform of ${p}(t,x)$ with respect to time $t$, written
in the travel time coordinates \eqref{eq:1D23},
\begin{equation}
\breve{p}(s,T) = \int_0^\infty  {p}(t,x(T)) e^{-s t} ds,
\end{equation}
satisfies the boundary value problem 
\begin{align}
\big(A_q + s^2) \breve{p}(s,T) &= s \sqrt{\sigma(T)} b(T) ~ ~ \mbox{for}~ T \in
(0, T_\ell),
\label{eq:SP1} \quad 
\partial_T \breve{p}(s,0) =  \breve{ p}(s,T_\ell) = 0.
\end{align}
Here we wrote the operator in \eqref{eq:1D5} 
in the travel time coordinates
\begin{equation}
A_q = - \sigma(T) \partial_T \Big( \frac{1}{\sigma(T)} \partial_T
\Big) = - \partial_T^2 + \partial_T q(T) \partial_T,
\label{eq:SP3}
\end{equation}
with $q(T) = \ln \sigma(T)$. This is the formulation considered in
\cite{borcea2005continuum}, and the spectral measure of $A_q$ is
defined by its eigenvalues $\{\lambda_j > 0\}_{j \ge 1}$ and
$\{\zeta_j = y_j^2(0)\}_{j \ge 1}$, where $\{y_j(T)\}_{j\ge 1}$ are
the eigenfunctions, orthonormal with respect to the weighted inner
product $\left< \cdot, \cdot \right>_{1/\sigma}$.

The spectrally matched ROM is defined in \cite{borcea2005continuum} by
an $n \times n$ tridiagonal matrix with spectral measure defined by
$\{\lambda_j, \zeta_j\}_{j=1}^n$.  To compare it with the ROM defined
in section \ref{sect:1D.2}, let us consider the Liouville transform
\begin{equation}
\breve{P}(s,T) = \frac{\breve{p}(s,T)}{\sqrt{\sigma(T)}},
\label{eq:SP4}
\end{equation}
and suppose that $\sigma(T)$ is constant  in a vicinity of $T = 0$. Then,
$\breve{P}(s,T) $ satisfies 
\begin{align}
\big(L_q L_q^T + s^2 \big)\breve{P}(s,T) &= s b(T) ~~ \mbox{for} ~ T
\in (0,T_\ell),
\label{eq:SP5} \quad 
\partial_T \breve{P}(s,0) = \breve{P}(s,T_\ell) = 0, 
\end{align}
with $L_q$ and $L_q^T$ defined in \eqref{eq:1D24} and
\eqref{eq:1D25}. Note that $L_q L_q^T$ is related to $A_q$ by a
similarity transformation,
\begin{equation}
L_q L_q^T = \sigma(T)^{-1/2} A_q \sigma(T)^{1/2},
\label{eq:SP6}
\end{equation}
so it has the same eigenvalues $\lambda_j$. The eigenfunctions $z_j(T)
= \sigma(T)^{-1/2} y_j(T)$ are orthonormal with respect to the
Euclidian $L^2\big((0,T_\ell)\big)$ inner product. Thus, the spectral
measure of $L_q L_q^T$ is the same as that of $A_q$, up to the multiplicative
constant $1/\sigma(0)$, which we take equal to $1$. 

The measurement
functions in the data model \eqref{eq:in5} are now
\begin{equation}
M_j(L_q L_q^T) = \left\{ \begin{array}{ll} \lambda_{j+1}, \quad & j =
  0, \ldots, n -1 \\ z^2_{j-n+1}(0), & j = n, \ldots 2n-1.
\end{array} \right.
\label{eq:SP7D}
\end{equation}
The ROM is defined by the $ n \times n $ symmetric, positive definite
and tridiagonal matrix $\tbA_q$ in the finite difference
discretization of \eqref{eq:SP5} on a special grid with $n$ points in
$[0,T_\ell]$,
\begin{equation}
\big(\tbA_q + s^2 \big)\tilde{\bp}(s) =  s \|b\| {\bf e}_1.\label{eq:SP7p}
\end{equation}
Let \begin{equation} \tbA_q = \tbL_q \tbL_q^T, \label{eq:CHSP} \end{equation} be the Cholesky factorization of
this matrix, with lower bidiagonal $\tbL_q$. Let also $\tilde
\lambda_j$ and $\tilde{\itbf z}_j$ be the eigenvectors of $\tbA_q$,
with Euclidean norm $\|\tilde{\itbf z}_j\|= 1$.  The matrix $ \tbA_q$
is obtained from the matching conditions
\begin{equation}
 M_j(L_q L_q^T) = \tilde M_j (\tbL_q \tbL_q^T) =
 \left\{ \begin{array}{ll} \tilde \lambda_{j+1}, \quad & j = 0,
   \ldots, n -1 \\ \big({\bf e}_1^T \tilde{\itbf z}_j\big)^2, & j = n,
   \ldots 2n-1,
\end{array} \right. 
\end{equation}
using the Lanczos algorithm \cite{chu2002structured}. The resulting
lower bidiagonal matrix $\tbL_q$ has the entries
\begin{align}
(\tbL_q)_{j,j} &= -\frac{1}{\sqrt{\gamma_j \hat \gamma_j}}, \quad 1 \le j \le n, \quad 
(\tbL_q)_{j+1,j} = \frac{1}{\sqrt{\gamma_{j} \hat
    \gamma_{j+1}}}, \quad 1 \le j \le n-1,
\label{eq:Sp9}
\end{align}
that have the same expression as in \eqref{eq:O3p}, but the values of
$\{\gamma_j, \hat \gamma_j\}_{j=1}^n$ are different.

\subsection{Inversion on the spectrally matched grid}
With the factorization \eqref{eq:CHSP} we can rewrite 
\eqref{eq:SP7p} as the first order system
\begin{align}
s \tilde \bcu(s) = \tbL_q^T \tilde \bp(s), \quad 
s \tilde \bp(s) = - \tbL_q \tilde \bcu(s) + \|b\|{\bf e}_1, 
\label{eq:Sp9p}
\end{align}
for the ROM primary and dual vectors $\tilde \bp(s)$ and $\tilde
\bcu(s)$.  We now show that the entries in these vectors represent
discretizations of the Laplace transforms $\breve{P}(s,T)$ and
$\breve{U}(s,T)$ of the primary and dual fields in equations
\eqref{eq:in.1N}, rewritten in travel time coordinates. The
discretization grid is defined by the spectrally matched ROM
coefficients $\{\gamma^0_j,\hat \gamma_j^0 \}_{j=1}^n$ calculated in
the reference medium with constant impedance $\sigma^0 = 1$ i.e.,
potential $q^0 = 0$.  It is proved in \cite[Lemma
  3.2]{borcea2005continuum} that these coefficients define a staggered
grid in the interval $[0,T_\ell]$,
\begin{equation}
0 = \hat T_0 = T_1  < \hat T_1 < T_2 < \hat T_2 < \ldots < \hat T_n < T_{n+1} ,
\label{eq:grid}
\end{equation}
with $T_{n+1}$ approaching $T_\ell$ from below in the limit $n \to
\infty$. The primary field $\breve P(s,T)$ is discretized on the grid
with points $\{T_j\}_{j=1}^{n+1}$, spaced at intervals $
h_j=\gamma_j^0 = T_{j+1}-T_j, $ and the dual field $\breve U(s,T)$ is
discretized on the grid with points $\{\hat T_j\}_{j=0}^{n}$, spaced
at intervals $ \hat h_j =\hat \gamma_j^0 = \hat T_{j}- \hat T_{j-1},
$ for $j = 1, \ldots, n$.

Let us define the diagonal matrices 
\[
H^{1/2}= \mbox{diag} \big(h_1^{1/2}, \ldots, h_n^{1/2}\big), \quad
\hat H^{1/2} = \mbox{diag} \big(\hat h_1^{1/2}, \ldots, \hat
h_n^{1/2}\big),
\] and write the vectors in \eqref{eq:Sp9p} in the form 
\begin{equation}
\tilde \bp(s) = \hat H^{1/2} \begin{pmatrix} \breve
  P_{_{T_1}}\hspace{-0.03in}(s) \\ \vdots \\ \breve
  P_{_{T_n}}\hspace{-0.03in}(s) \end{pmatrix}, \quad \tilde \bcu(s) =
H^{1/2} \begin{pmatrix} \breve U_{_{\hat T_1}}\hspace{-0.03in}(s)
  \\ \vdots \\ \breve U_{_{\hat
      T_n}}\hspace{-0.03in}(s) \end{pmatrix},
\end{equation}
so that the first equation in the system  \eqref{eq:Sp9p} becomes
\begin{align}
s \begin{pmatrix} \breve U_{_{\hat T_1}}\hspace{-0.03in}(s) \\ \vdots \\ 
\breve U_{_{\hat T_n}}\hspace{-0.03in}(s) \end{pmatrix} &= H^{-1/2} \tbL_q^T \hat H^{1/2} 
\begin{pmatrix} \breve P_{_{T_1}}\hspace{-0.03in}(s) \\ \vdots \\ 
\breve P_{_{T_n}}\hspace{-0.03in}(s)
\end{pmatrix}. \label{eq:SP10p}
\end{align}
Let also
\begin{equation}
\sigma_j = \frac{\hat h_j}{\hat \gamma_j}= \frac{\hat \gamma_j^0}{\hat
  \gamma_j}, \quad \hat \sigma_j = \frac{\gamma_j}{h_j} =
\frac{\gamma_j}{\gamma_j^0}, \quad j = 1, \ldots, n,
\label{eq:ratios}
\end{equation}
and write explicitly  the $j$-th equation in \eqref{eq:SP10p}
\begin{align*}
s \, \breve U_{_{\hat T_j}}\hspace{-0.03in}(s) = \frac{\breve
  P_{_{T_{j+1}}}\hspace{-0.03in}(s)- \breve
  P_{_{T_j}}\hspace{-0.03in}(s)}{h_j} + \breve
P_{_{T_{j+1}}}\hspace{-0.03in}(s)
\left(\frac{\sqrt{\sigma_{j+1}}-\sqrt{\hat \sigma_j}}{h_j \sqrt{\hat
    \sigma_j}} \right) + \breve P_{_{T_j}}\hspace{-0.03in}(s)
\left(\frac{\sqrt{\hat \sigma_{j}}-\sqrt{ \sigma_j}}{h_j \sqrt{\hat
    \sigma_j}}\right).
\end{align*}
This is the discretization of equation
\begin{equation*}
s \breve U(s,T) = L_q^T \breve P(s,T) = \partial_T \breve P(s,T) +
\partial_T \ln \sqrt{\sigma(T)} \breve P(s,T),
\end{equation*}
on the spectrally matched grid. A similar result applies to the second
equation in \eqref{eq:Sp9p}, which is the discretization of
\begin{equation*}
s \breve P(s,T) = -L_q \breve U(s,T) + b(T).
\end{equation*}

The ratios \eqref{eq:ratios} are of the same form as in
\eqref{eq:ratiosROM} and represent approximations of the impedance
function $\sigma(T)$ on the spectrally matched grid
\eqref{eq:grid}. Specifically, it is proved in \cite[Theorem
  6.1]{borcea2005continuum} that if $\sigma(T)$ is smooth, then the
piecewise constant function
\[
\sigma^n(T) = \left\{ \begin{array}{ll} \sigma_j \quad & \mbox{for} ~
  T \in [T_j,\hat T_j), ~ j = 1, \ldots, n, \\ \hat \sigma_j &
    \mbox{for} ~ T \in [\hat T_j, T_{j+1}), ~ j = 1, \ldots, n,
      \\ \hat \sigma_n & \mbox{for} ~ T \in [T_{n+1},T_\ell],
\end{array} \right.
\]
converges to $\sigma(T)$ pointwise and in $L^1[0,T_\ell]$, as $n \to \infty$.

\section{The multi dimensional case}
\label{sect:MD}
In this section we generalize the DtB mapping from one dimension, 
as described in section \ref{sect:1D}, to $\mathbb{R}^d$ with  $d > 1$.
The derivation follows the same strategy as in section \ref{sect:1D}, with  
certain modifications described below.

\subsection{Data model for an array of sensors}
\label{sect:dataarray}

In the multi-dimensional case we consider an array of $m$ sensors on
the accessible boundary $\partial \Omega_a$, located at points $\bx_s
\in \partial \Omega_a$, for $s = 1,\ldots,m$. Each sensor excites a
pressure field, the solution of the wave equation
\begin{equation}
\big( \partial_t^2 + A \big) \mathfrak{p}_{s}(t, \bx) = \partial_t
f(t) \delta(\bx-\bx_s), \quad - \infty < t < \infty, ~ ~ \bx \in
\Omega,
\label{eq:multiwaveps}
\end{equation}
where the index $s$ denotes the source and the operator $A$ is now
defined as
\begin{equation}
A = - \sigma(\bx) c(\bx) \nabla \cdot \Big[ \frac{c(\bx)}{\sigma(\bx)}
  \nabla \Big],
\label{eq:multiA}
\end{equation}
with $\nabla$ denoting the gradient and $\nabla \cdot$ the divergence operator.
For simplicity, we assume that the same pulse $\partial_t f(t)$ is
emitted from all the sensors.
The boundary conditions at $\partial \Omega_a$ are $ \nabla
\mathfrak{p}_s(t,\bx) \cdot \bn(\bx) = 0, $ where $\bn(\bx)$ is the
outer unit normal, and on the inaccessible boundary $\partial
\Omega_i$ we let $ \mathfrak{p}_s(t,\bx) = 0.  $ The medium is at rest
initially, so we set $\mathfrak{p}_s(t, \bx) = 0,$ for $ t \ll 0.$

Following the same argument that lead to equation \eqref{eq:1D11}, we
define the matrix-valued data $\bD_k \in \mathbb{R}^{m \times m}$ with entries 
defined by 
\begin{equation}
(\bD_k)_{i,j} = \mathfrak{p}_i(t_k, \bx_j) + \mathfrak{p}_i(-t_k,
  \bx_j), \quad i,j = 1,\ldots,m,
\end{equation}
in terms of the measurements at the  instances  $t_k = k \tau$ for $k = 0, 1, \ldots, 2n-1$.   For each $k$ the matrix $\bD_k$ is symmetric
due to the source-receiver reciprocity.

\subsection{First order system form and Liouville transformation}
\label{sect:multifo}

Similar to the one-dimensional case, we introduce the sensor functions
\begin{equation}
b_s(\bx) = \sqrt{\frac{\sigma(\bx_s) c(\bx_s)}{\sigma(\bx)}} 
[\hat f(\sqrt{A})]^{1/2} \delta(\bx - \bx_s), 
\quad \bx \in \Omega, \quad s = 1,\ldots,m,
\end{equation}
and define the analogue of \eqref{eq:1D14}
\begin{equation}
{p}_s(t,\bx) = \cos(t \sqrt{A}) \sqrt{\sigma(\bx)}
b_s(\bx), \quad s=1,\ldots,m.
\label{eq:multipfrakcos}
\end{equation}
Here $p_s(t,\bx)$ is the pressure field in the first order system
\begin{align}
\partial_t \begin{pmatrix} {p}_s(t,\bx) \\ -{\itbf u}_s(t, \bx) \end{pmatrix} = 
\begin{pmatrix} 0 & \sigma(\bx) c(\bx) \nabla \cdot \\
\frac{c(\bx)}{\sigma(\bx)} \nabla & 0 \end{pmatrix} \begin{pmatrix}
  {p}_s(t,\bx) \\-{\itbf u}_s(t,\bx) \end{pmatrix},
\quad t > 0, \quad \bx \in \Omega,
\label{eq:multifo1}
\end{align}
with initial conditions 
\begin{equation}
{p}_s(0,\bx) = \sqrt{\sigma(\bx)}
\cb_s(\bx) \quad \mbox{and} \quad   {\itbf u}_s(0,\bx) = 0, \qquad  \bx \in \Omega,
\end{equation}
and
boundary conditions  
\begin{equation}
{p}_s(t,\bx)|_{\bx \in \Omega_i} = 0 \quad \mbox{and} \quad 
{\itbf u}_s(t,\bx)|_{\bx \in \Omega_a} = 0, \qquad t > 0.
\end{equation}  
This system
is the analogue of \eqref{eq:1D15} and the vector field ${\itbf
  u}_s(t,\bx)$ is the particle velocity.

Using the Liouville transformation 
\begin{equation}
\cp(t,\bx;\bx_s) = \frac{p(t,\bx;\bx_s)}{\sqrt{\sigma(\bx)}}, \quad 
\bcu(t,\bx;\bx_s) = -\sqrt{\sigma(\bx)} {\itbf u}(t,\bx;\bx_s),
\label{eq:multicpcu}
\end{equation}
we rewrite \eqref{eq:multifo1} as a Schr\"{o}dinger system
\eqref{eq:in.1} with the operators $L_q$ and $L_q^T$ given by
\begin{align}
L_q  = - c(\bx) \nabla \cdot \; + \frac{1}{2} c(\bx) \nabla q(\bx) \cdot
 ~ ~ \mbox{and} ~ ~ 
\label{eq:multilq}
L_q^T & = c(\bx) \nabla + \frac{1}{2} c(\bx) \nabla q(\bx),
\end{align}
and the same  $q(\bx) = \ln \sigma (\bx)$.
The transformed boundary and initial conditions take the form
\eqref{eq:in.3}--\eqref{eq:in.2}, and the entries of the data matrix
$\bD_k \in \mathbb{R}^{m \times m}$, for $k = 0, \ldots, 2n-1$, are
expressed in terms of the primary wave $\cp_i(t,\bx)$ as
\begin{equation}
(\bD_k)_{i,j} = \left< b_j, \cp_i(t_k,\cdot) \right>_{\frac{1}{c(\bx)}} = 
\int_\Omega \frac{d\bx}{c(\bx)} \cb_j(\bx) \cp_i(t_k, \bx), 
\quad i,j = 1,\ldots,m.
\label{eq:multidatab}
\end{equation}

\subsection{Symmetrized data model, propagator and measurement function}
\label{sect:multisym}
In one dimension we used travel time coordinates to write the data
model in the symmetrized form. Because such a transformation is not
available in higher dimensions, we follow a different approach to
symmetrize \eqref{eq:multidatab} and thus obtain an analogue of
\eqref{eq:1D33}.

Combining \eqref{eq:multidatab} with \eqref{eq:multipfrakcos} and
\eqref{eq:multicpcu} we write
\begin{equation}
(\bD_k)_{i,j} = \int_{\Omega} \frac{d \bx}{c(\bx)} 
b_j(\bx) \sigma^{-1/2}(\bx) \cos(t_k \sqrt{A}) \sigma^{1/2}(\bx) b_i(\bx),
\end{equation}
where we assume for the remainder of this section $i,j = 1,\ldots,m$,
and $k=0,\ldots,2n-1$.  From the definition \eqref{eq:multilq} of
$L_q$ it follows that \begin{equation} A = \sigma^{1/2}(\bx) L_q L_q^T
\sigma^{-1/2}(\bx).\end{equation} Hence,
\begin{equation}
(\bD_k)_{i,j} = \int_{\Omega} \frac{d \bx}{c(\bx)} b_j(\bx) 
\cos \Big( t_k \sqrt{L_q L_q^T} \Big) b_i(\bx),
\label{eq:multidatac}
\end{equation}
where we used that analytic matrix functions commute with similarity
transformations.  We use another similarity transformation to
rewrite \eqref{eq:multidatac} as
\begin{equation}
(\bD_k)_{i,j} = \int_{\Omega} d \bx \; b_j^c(\bx) 
\cos\Big(t_k \sqrt{c^{-1/2} L_q L_q^T c^{1/2}} \Big) b_i^c(\bx),
\label{eq:multidatad}
\end{equation}
with the rescaled sensor functions 
\begin{equation}
b^c_s (\bx) = c^{-1/2}(\bx) b_s(\bx), \quad s = 1,\ldots,m.
\end{equation}
We also define the rescaled operators
\begin{align}
L_q^c & = c^{-1/2}(\bx) L_q c^{1/2}(\bx)= c^{1/2}(\bx) \left( -
\nabla \cdot + \frac{1}{2} \nabla q(\bx) \; \cdot \right)
c^{1/2}(\bx),
\label{eq:multilqsym} \\
L_q^{c^T} & = c^{-1/2}(\bx) L_q^T c^{1/2}(\bx) = c^{1/2}(\bx) \left(
\nabla + \frac{1}{2} \nabla q(\bx) \right) c^{1/2}(\bx),
\label{eq:multilqtsym}
\end{align}
which are adjoint to each other with respect to the standard
$L^2(\Omega)$ inner product.  These operators define the propagator
for the multi-dimensional problem as in \eqref{eq:1D29},
\begin{equation}
\msP = \cos \Big( \tau \sqrt{L_q^c L_q^{c^T}} \Big).
\label{eq:multiprop}
\end{equation}
The data model \eqref{eq:multidatad} is now in symmetric form, and the
measurement functions $\bM_k(L_q L_q^T) \in \mathbb{R}^{m
  \times m}$ are defined component-wise by 
\begin{equation}
(\bD_k)_{i,j} = [\bM_k (L_q^c L_q^{c^T})]_{i,j} = 
\left< b_j^c(\bx), \mathcal{T}_k(\msP) b_i^c(\bx) \right>.
\label{eq:DataMulti}
\end{equation}

Similar to the one-dimensional case, the propagator $\msP$ can be used
to define an exact time stepping scheme for
\begin{equation}
\label{eq:multisnap}
\cp_{k,s}(\bx) = \mathcal{T}_k (\msP) \cb_s^c(\bx), \quad
k=0,\ldots,2n-1, \quad s=1,\ldots,m.
\end{equation}
Here we use the convention that the first index  denotes
the discrete time instance and the second index denotes the
source.  We obtain 
the same second-order time stepping scheme \eqref{eq:1D39},
\begin{equation}
\frac{1}{\tau^2}\Big[\cp_{k+1,s} (\bx) - 2 \cp_{k,s} (\bx) + \cp_{
    k-1,s}(\bx)\Big] = -\xi(\msP) \cp_{k,s}(\bx), \quad
k=0,\ldots,2n-1,
\label{eq:multitimestep}
\end{equation}
with the affine function $\xi$ given by \eqref{eq:defXi} and the
initial conditions 
\begin{equation}
 \cp_{0,s} (\bx) = \cb^c_s(\bx), \quad \cp_{-1,s}
(\bx) = \cp_{1,s} (\bx) = \msP \cb^c_s (\bx).
\end{equation} 
We also have the same
factorization \eqref{eq:1D41} of $\xi(\msP)$ in terms of the operator
$\mathcal{L}_q$ defined as in \eqref{eq:1D42}, with $L_q$ and $L_q^T$
replaced by $L_q^c$ and $L_q^{c^T}$.

\subsection{Multi-input, multi-output reduced order model}
\label{sect:mimorom}

The main difference between one and multi dimensions is the type of
ROM that we use. In one dimension we had a single-input, single output
(SISO) projection ROM \eqref{eq:1D56p}, obtained from matching the ROM
output \eqref{eq:ROMMeas} to the scalar valued data \eqref{eq:1D50}.
In multi dimensions we need a multi-input, multi-output (MIMO) ROM that 
matches the matrix valued data \eqref{eq:DataMulti}.

As in section \ref{sect:1D.2}, let us introduce the matrix $\bA \in
\mathbb{R}^{N \times N}$, a discretization of the operator
\eqref{eq:multiA} on a very fine, uniform grid with a total of $N$
nodes and step size $h$. Note that the operator \eqref{eq:multiA} is
related to the Schr\"{o}dinger operators
\eqref{eq:multilqsym}--\eqref{eq:multilqtsym} as
\begin{equation}
A = \sigma^{1/2}(\bx) L_q L_q^T \sigma^{-1/2}(\bx) = [\sigma(\bx)
  c(\bx)]^{1/2} L_q^c L_q^{c^T} [\sigma(\bx) c(\bx)]^{-1/2}.
\end{equation}
Let $\boldsymbol \Sigma \in \mathbb{R}^{N \times N}$ and $\bC \in
\mathbb{R}^{N \times N}$ be the diagonal matrices with entries given
by $\sigma(\bx)$ and $c(\bx)$ evaluated at the fine grid nodes. Then
we can set $\bL_q$, a fine grid approximation of $L_q^c$, to be a
Cholesky factor of
\begin{equation}
{\boldsymbol \Sigma}^{-1/2} \bC^{-1/2} \bA {\boldsymbol \Sigma}^{1/2} \bC^{1/2} 
=  \bL_q \bL_q^T,
\end{equation}
where we drop the index $c$ on $\bL_q$ to simplify notation.  We
assume as in one dimension that the discretized propagator $\bbP \in
\mathbb{R}^{N \times N}$ is well approximated by \eqref{eq:approxP} on
the fine grid.

As in one dimension, we call the fine grid discretization of the field $P$ at the measurement 
instances "the primary snapshots". It is convenient to arrange these into matrices
$\bp_k \in \mathbb{R}^{N \times m}$, for $k = 0,\ldots,2n-1$, with
each column corresponding to a different sensor. These matrices
satisfy a fine grid analogue of the time stepping scheme
\eqref{eq:multitimestep},
\begin{equation}
\frac{1}{\tau^2}\big(\bp_{k+1} - 2 \bp_k + \bp_{k-1}\big) = 
-\xi(\bbP) \bp_k, \quad k=0,\ldots,2n-1,
\label{eq:multitimesteppmat}
\end{equation}
with initial conditions $ \bp_0 = \bb$ and $\bp_{-1} = \bp_1 = \bbP
\bb.  $ The sensor matrix $\bb \in \mathbb{R}^{N \times m}$ is 
\begin{equation}
\bb = [\bb_1,\ldots,\bb_m],
\label{eqn:bb}
\end{equation}
where the entries of each column $\bb_s \in \mathbb{R}^N$ are the
values of the rescaled sensor function $b_s^c(\bx)$ evaluated on the
fine grid, multiplied by $h^{d/2}$.

The discretized data model and the measurement functions are given by 
\begin{equation}
\bD_k = \bM_k(L_q L_q^T) \approx \bM_k (\bL_q \bL_q^T) = \bb^T \bp_k =
\bb^T \mathcal{T}_k (\bbP) \bb \in \mathbb{R}^{m \times m},
\label{eq:multidatadisc}
\end{equation}
for $k=0,\ldots,2n-1$. The analogue of equation \eqref{eq:bb}, which
relates the data $\bD_0$ at the first time instant to the sensor
matrix is 
\begin{equation}
\bD_0 = 
\begin{bmatrix} 
\left< b^c_1(\bx), b^c_1(\bx) \right> & \left< b^c_1(\bx), b^c_2(\bx) \right> & 
\ldots & \left< b^c_1(\bx), b^c_m(\bx) \right> \\ 
\vdots & \vdots & \ddots & \vdots \\
\left< b^c_m(\bx), b^c_1(\bx) \right> & \left< b^c_m(\bx), b^c_2(\bx) \right> & 
\ldots & \left< b^c_m(\bx), b^c_m(\bx) \right> 
\end{bmatrix} \approx 
\bb^T \bb.
\label{eq:multibb}
\end{equation}
As we did in one dimension, we neglect henceforth the fine grid
discretization errors and treat the approximate relations
\eqref{eq:multidatadisc}--\eqref{eq:multibb} as equalities.

The MIMO ROM consists of the symmetric matrix $\tilde{\bbP} \in
\mathbb{R}^{nm \times nm}$, called the ROM propagator, and the ROM
sensor matrix $\tilde \bb \in \mathbb{R}^{nm \times m}$, satisfying 
the data matching conditions
\begin{equation}
\bD_k = \tilde \bb^T \mathcal{T}_k (\tilde{\bbP} ) \tilde \bb, 
\quad k = 0, \ldots, 2n-1.
\label{eq:multidatamatch}
\end{equation}
The matrix $\tilde{\bbP}$ is block tridiagonal, with $m \times m$
blocks, while $\tilde \bb$ has all zeros
except for the uppermost $m \times m$ block. Using the
$m\times m$ identity matrix $\bI_m$ and zero matrix $\boldsymbol{0}_m$, 
we write 
\begin{equation}
\tilde \bb =  \bE_1 (\bD_0)^{1/2},
~ ~ \mbox{with} ~ ~
\bE_1 = \begin{bmatrix} 
\bI_m \\ {\boldsymbol 0}_m \\ \vdots \\ {\boldsymbol 0}_m 
\end{bmatrix} \in \mathbb{R}^{nm \times m}.
\label{eq:e1mat}
\end{equation}

\subsection{Calculation of the projection MIMO ROM}
\label{sect:mimoromcalc}
The MIMO ROM satisfying the data matching conditions
\eqref{eq:multidatamatch} is given by the orthogonal projection of
$\bbP$ on the block Krylov subspace, spanned by the columns of \[ \bP =
      [\bp_0, \ldots, \bp_{n-1}] \in \mathbb{R}^{N \times nm}.\] We
      follow the notation in section \ref{sect:projROM}, and let $\bV
      = [\bV_1,\ldots,\bV_n]$ be the matrix containing the orthonormal
      basis for $\mbox{range}(\bP)$. Here $\bV \in \mathbb{R}^{N
        \times nm}$ and each $\bV_k$ is an $N \times m$ matrix.

To compute the matrices $\bP^T \bP$ and $\bP^T \bbP \bP$ from the data
we can still use the formulas \eqref{eq:MassM1} and \eqref{eq:Stiff2},
however, the indexing is understood block-wise. Thus,
when we write
\begin{align}
\big(\bP^T \bP\big)_{i,j} & = \frac{1}{2} \Big( \bD_{i+j-2} +
\bD_{|i-j|} \Big) \in \mathbb{R}^{m \times
  m}, \label{eq:multimassdata} \\ (\bP^T \bbP \bP)_{i,j} & =
\frac{1}{4} \Big( \bD_{i+j-1} + \bD_{|j-i+1|} + \bD_{|j-i-1|} +
\bD_{|j+i-3|} \Big) \in \mathbb{R}^{m \times
  m}, \label{eq:multistiffdata}
\end{align}
for $i,j = 1,\ldots,n$, we use the notation $\big( \bP^T \bP
\big)_{i,j}$ for the $m \times m$ block of $\bP^T \bP \in
\mathbb{R}^{nm \times nm}$, at the intersection of rows
$(i-1)m+1,\ldots,im$ and columns $(j-1)m+1,\ldots,jm$. We use 
this notation for $\bP^T \bbP \bP \in \mathbb{R}^{nm \times
  nm}$ and all other block matrices with block size $m \times m$. 

Note that \eqref{eq:multimassdata}--\eqref{eq:multistiffdata} give
that $\bP^T \bP$ and $\bP^T \bbP \bP$ and their blocks are symmetric.

\subsection{The Data to Born map} 
\label{sect:mimoromfactor}
The main difference in the calculation of the DtB mapping is that the
Cholesky factorizations \eqref{eq:Chol1} and \eqref{eq:Chol} at step 1
of Algorithm \ref{eq:algo1} are replaced with their block Cholesky
counterparts given below.

\vspace{0.05in}\begin{alg}[Block Cholesky factorization]
\label{alg:blockhol}~\\
\textbf{Input:} the symmetric block matrix $\bX \in \mathbb{R}^{nm \times
  nm}$ with $m \times m$ blocks.

\vspace{0.05in} \noindent To obtain the block Cholesky factorization
of $\bX$ perform the following steps:
\begin{equation}
\hspace{-1in} \mbox{For} ~ k = 0, 1, \ldots, n-1 ~ \mbox{compute:} ~ \bR_{k,k} =
\bQ_k\Big( \bX_{k,k} - \sum_{i=0}^{k-1} \bR_{i,k}^T \bR_{i,k}
\Big)^{1/2}, \label{eq:diagblock}
\end{equation}
where $\bQ_k\in \mathbb{R}^{m\times m}$ is an arbitrary orthogonal
matrix. 
\begin{equation}
\hspace{-1in} \mbox{For} ~ j = k+1, \ldots, n-1 ~ \mbox{compute:} ~
\bR_{k,j} = \bR_{k,k}^{-1} \Big( \bX_{k,j} - \sum_{i=0}^{k-1}
\bR_{i,k}^T \bR_{i,j} \Big).
\end{equation}
 \textbf{Output:} the block matrix $\bR \in
\mathbb{R}^{nm \times nm}$ with $m \times m$ blocks, satisfying $\bX =
\bR^T \bR$.
\end{alg}

\vspace{0.05in} While the regular Cholesky factorization is defined
uniquely (assuming it uses the principal value of the square root),
there is an ambiguity in defining the block Cholesky factorization,
which comes from the computation of the diagonal blocks $\bR_{k,k}$ in
\eqref{eq:diagblock}. An optimal choice of this matrix is still an
open question. We obtained good results with $\bQ_k = {\itbf I}$, but
here we present another choice (yielding equally good numerical
results). This choice allows us to extend the Galerkin-Petrov
reasoning of Section~\ref{sect:1D.3} to the MIMO case.  Explicitly, we
choose the factor $\bQ_k$ consistent with the MIMO analogue of
recursion \eqref{eq:O1}--\eqref{eq:O3p} for computing the primary and
dual orthogonalized block snapshots $\overline{\bf p}_{j}\in
\mathbb{R}^{N \times m}$ and $\overline{\bf u}_{j}\in \mathbb{R}^{N
  \times m}$, for $j = 1,\ldots,n.$ Here the coefficients are no
longer scalar, but symmetric positive definite matrices $\bgamma_j \in
\mathbb{R}^{m \times m}$ and $\bgammahat_j \in \mathbb{R}^{m \times
  m}$.  These matrices give the estimates of $\sigma(\bx)$ in the
multidimensional case \cite[Section 7.3]{druskin2016direct}, so we can
extend the reasoning of section~\ref{sect:1DSPECT} to the MIMO case,
relating the block-bidiagonal $\tbL_q$ to a discretization of the
Schrodinger equation, via the solution of the discrete inverse
problem. We describe the computation of the block-bidiagonal $\tbL_q$
in Appendix~\ref{ap:blockCholesky}.

After the factor $\tbL_q$ is found, the DtB map is computed using the
block versions of Algorithms \ref{eq:algo1} and \ref{alg:chain}, with
$\bb$ given by \eqref{eqn:bb} and ${\bf e}_1$ replaced by $\bE_1$ from
\eqref{eq:e1mat}.  We can also rewrite the MIMO counterpart of
\eqref{eq:GP2} in block form. The validity of the MIMO DtB map is
based on the assumption of weak dependence of the primary and dual
block-QR orthonormal snapshots $\bV$ and $\bW$ on $q = \ln \sigma$.
Because of the consistency with the discrete inverse problem discussed
above, this weak dependence can be understood using the same reasoning
as in the SISO case. In particular, similar to the SISO case, the
orthonormal snapshots approximate columns of the identity, as shown in
Figure~\ref{fig:UV}. We also refer for more details to
\cite{DruskinMamonovZaslavsky2017}.

\begin{rem}
Unlike in the one dimensional case, in multi
dimensions the DtB algorithm becomes mildly ill-posed, even for the
space-time sampling close to the Nyquist rate. A simple regularization
algorithm presented in \cite{DruskinMamonovZaslavsky2017} makes the
DtB mapping practically insensitive to a reasonable (order of few per
cent) level of noise in the measured data for the problem sizes
considered in the numerical simulations shown in the next section.
\end{rem}

\subsection{Numerical results}
\label{sect:num2D}
\begin{figure}[t!]
\centering
\begin{tabular}{cc}
$\sigma(\bx)$ & $c(\bx)$ \\
\includegraphics[width=0.45\textwidth]{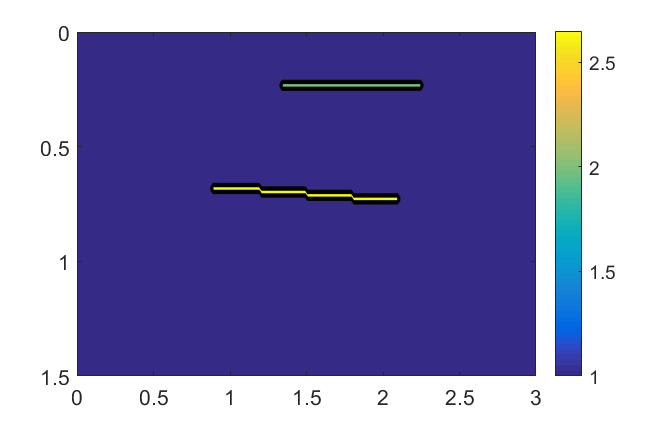} & 
\includegraphics[width=0.445\textwidth]{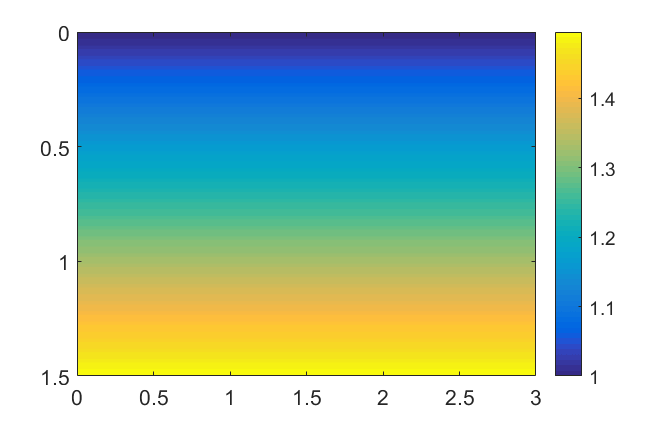}
\end{tabular}
\vspace{-0.12in}\caption{Model $\sigma(\bx)$ with two inclusions (left) and linear
  $c(\bx)$ (right). The axes are in km units.}
\label{fig:model2d}
\end{figure}

We begin with numerical results for a two dimensional impedance model
with two inclusions and a linear velocity model, shown in
Figure~\ref{fig:model2d}. We display the relative impedance and wave
speed, normalized by their constant values at the sensors. Both the
time sampling $\tau$ and the distances between the $m=50$ sensors in
the array are chosen close to the Nyquist sampling rate for the
Gaussian pulse used in the experiments. As in the one dimensional
case, the scattering data and the true Born approximation are computed
using a fine grid finite difference time domain scheme, with grid
steps much smaller than $\tau$.

\begin{figure}
\begin{tabular}{cccc}
Columns of $\bP^0$ & Columns of $\bP$ & Columns of $\bV^0$ & Columns
of $\bV$
\\ \includegraphics[width=0.22\textwidth]{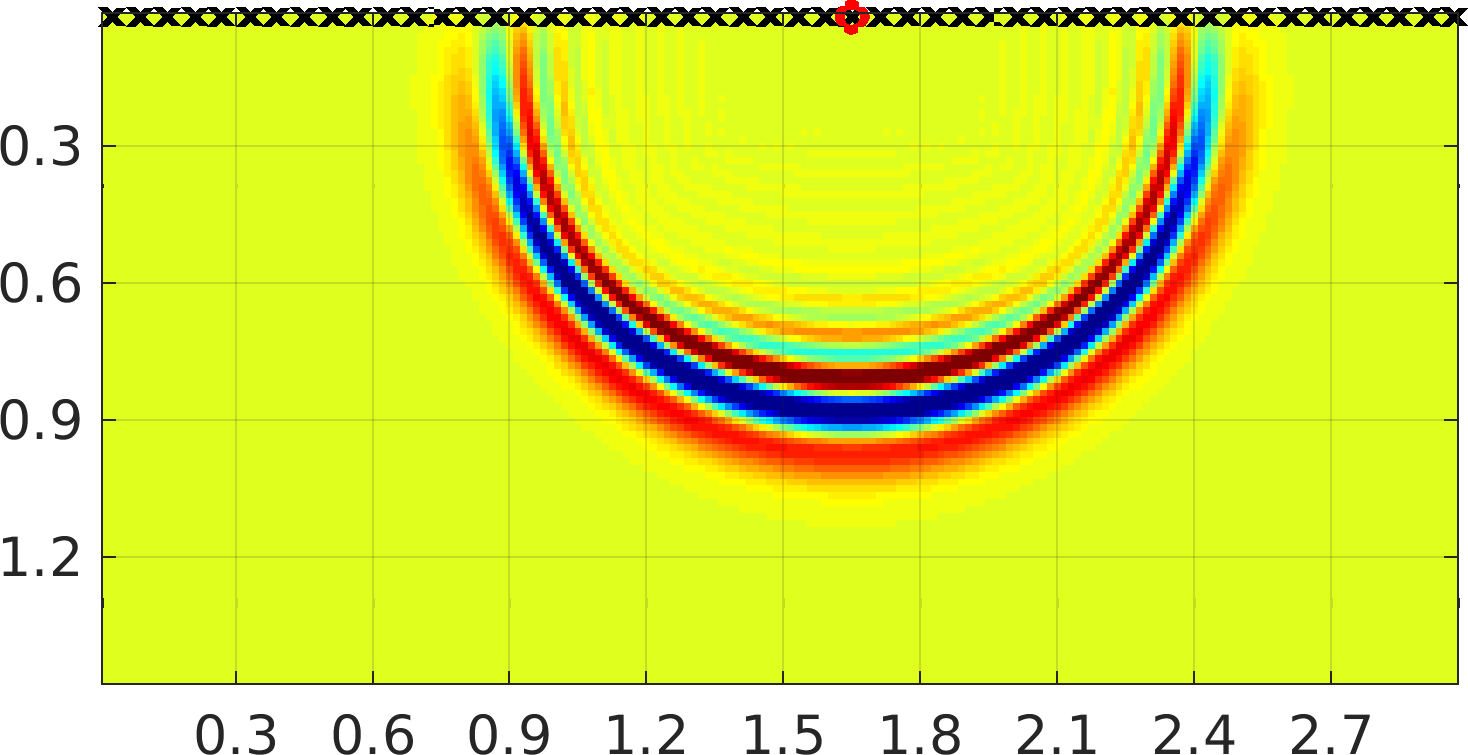}
& \includegraphics[width=0.22\textwidth]{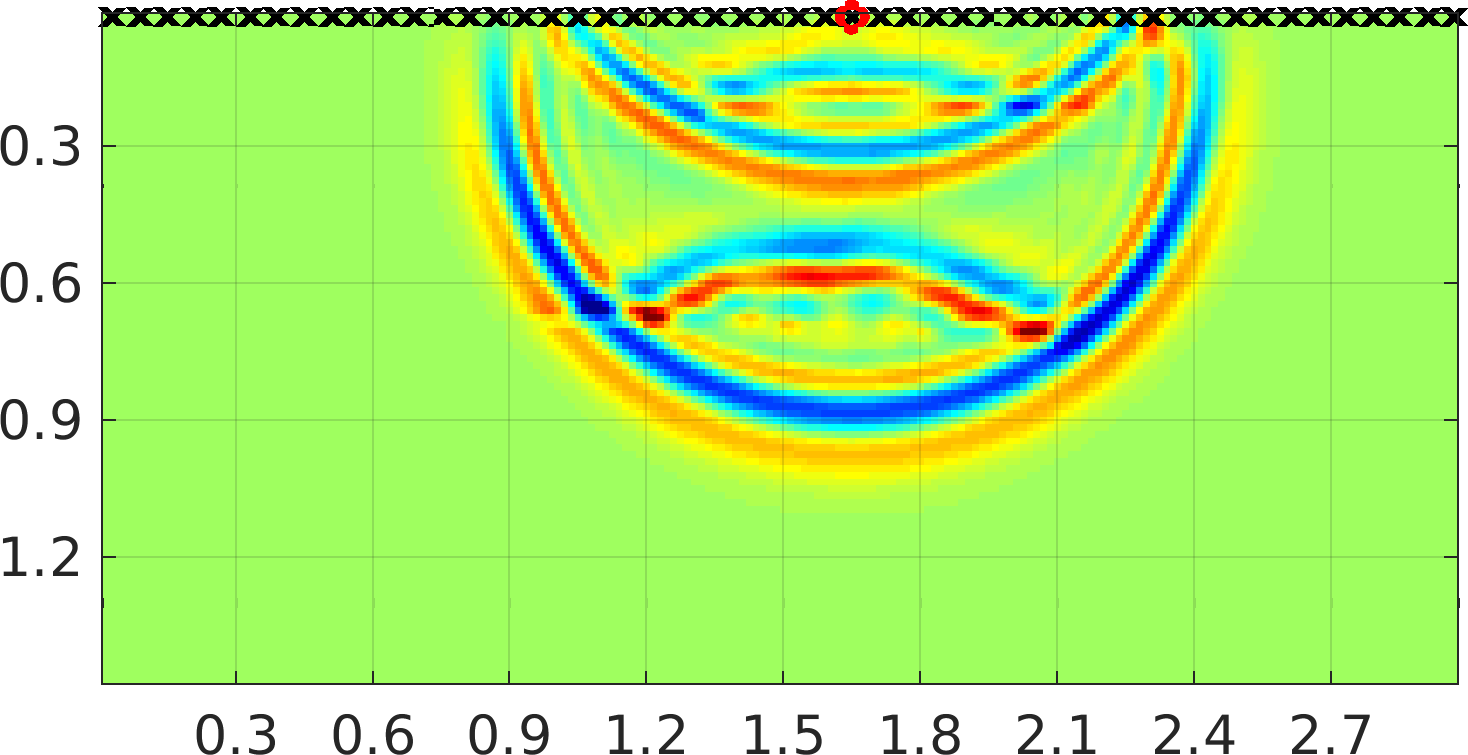} &
\includegraphics[width=0.22\textwidth]{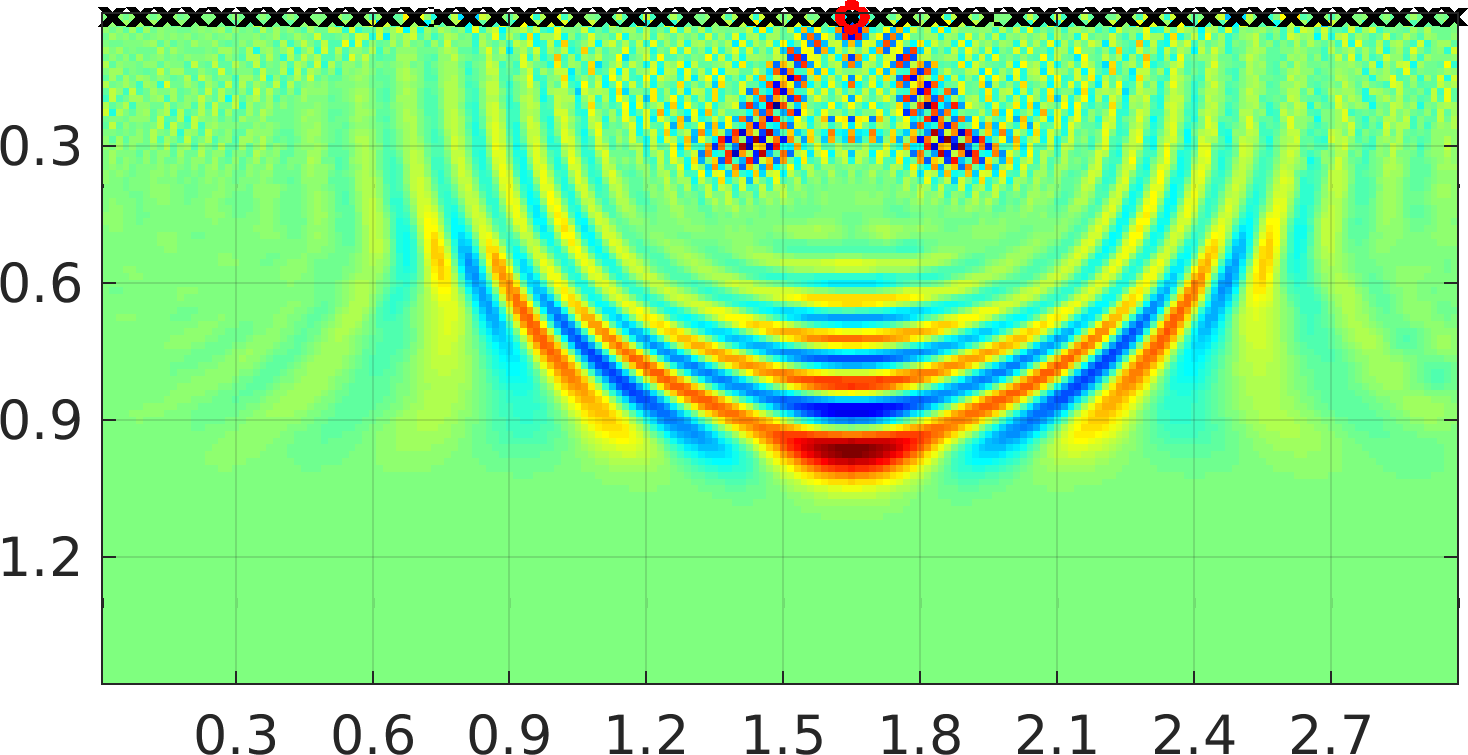} &
\includegraphics[width=0.22\textwidth]{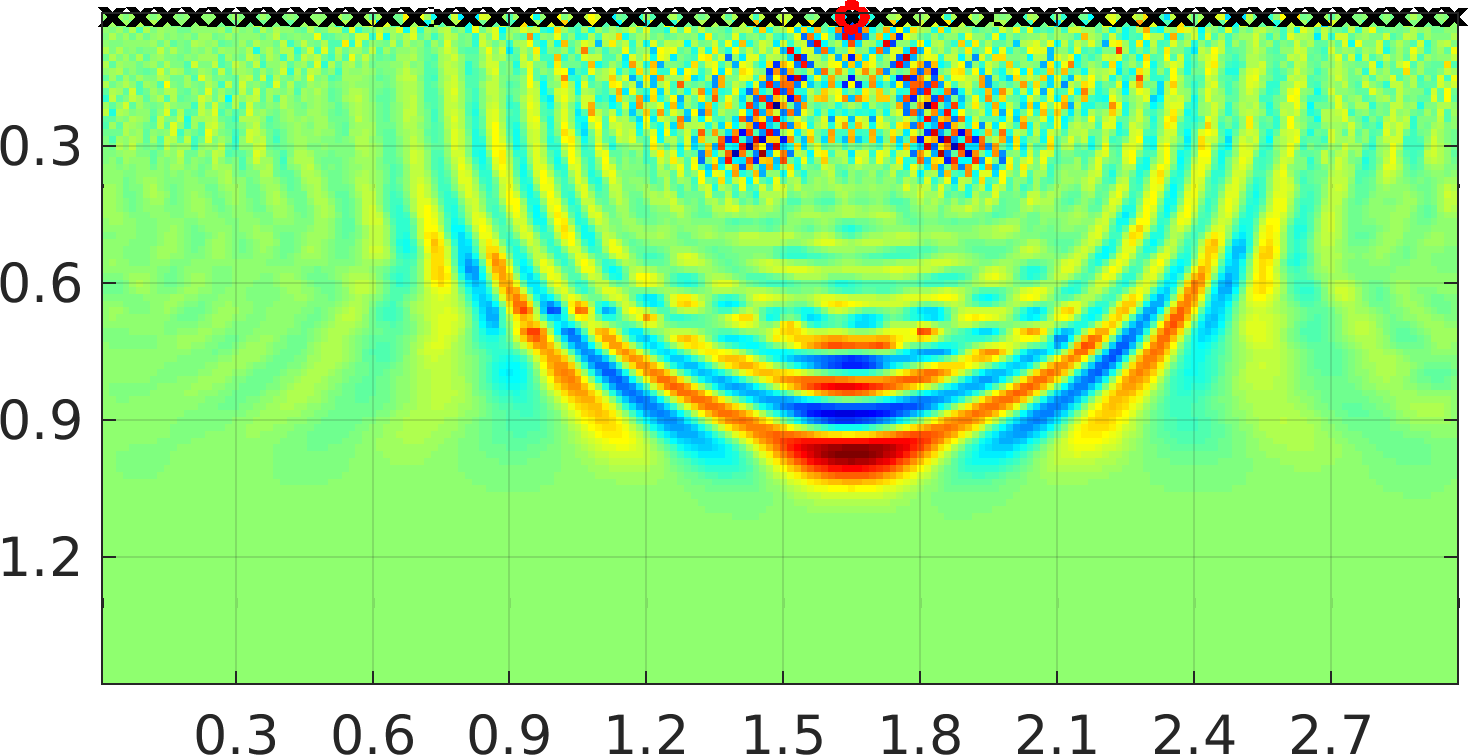}
\\ \includegraphics[width=0.22\textwidth]{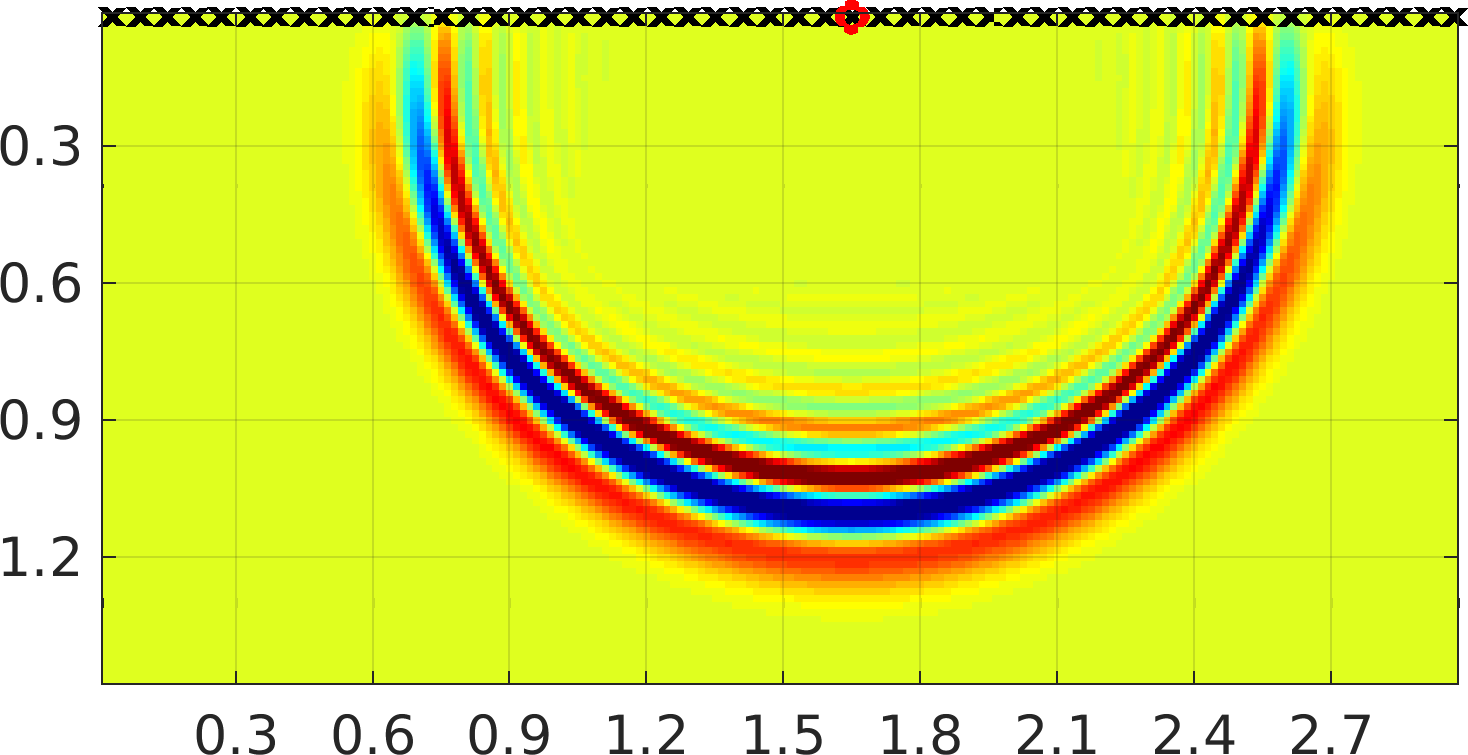}
& \includegraphics[width=0.22\textwidth]{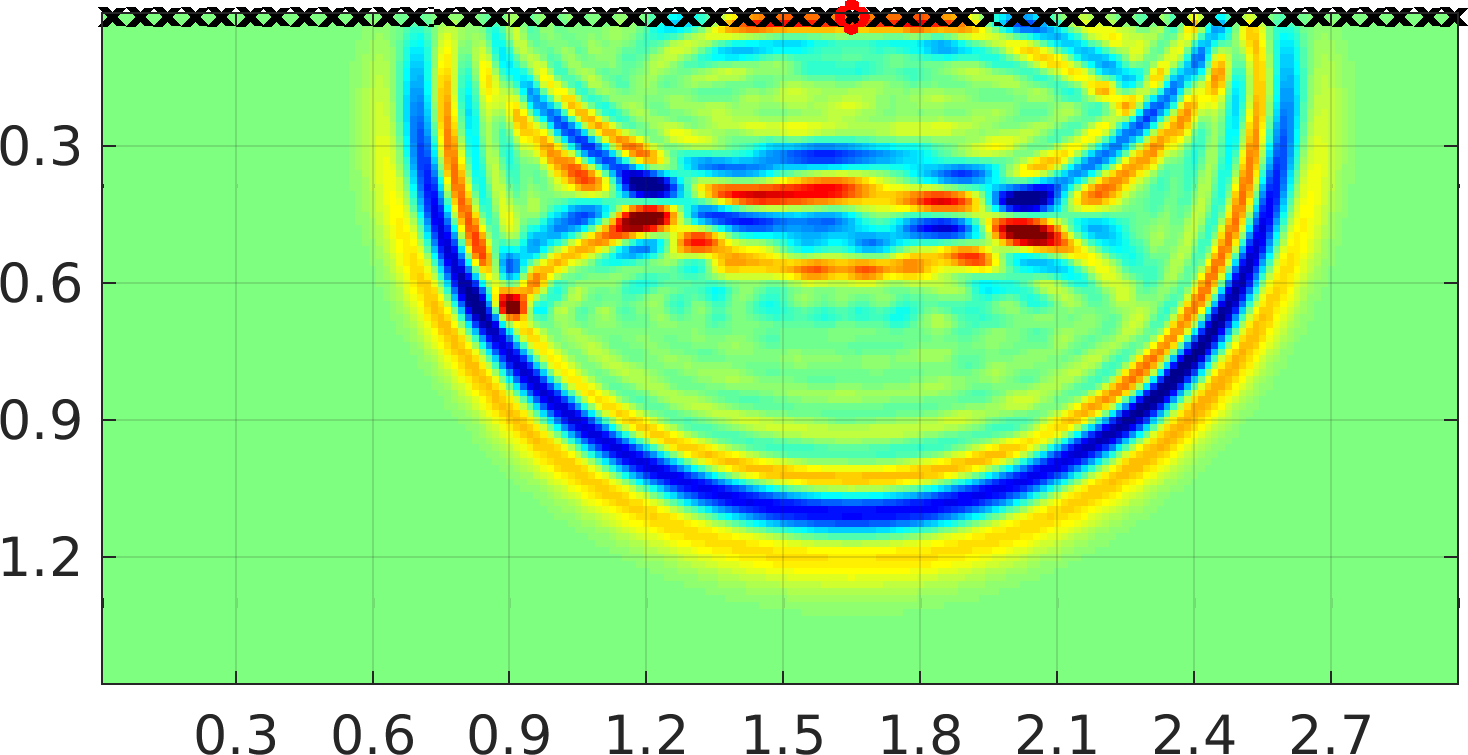} &
\includegraphics[width=0.22\textwidth]{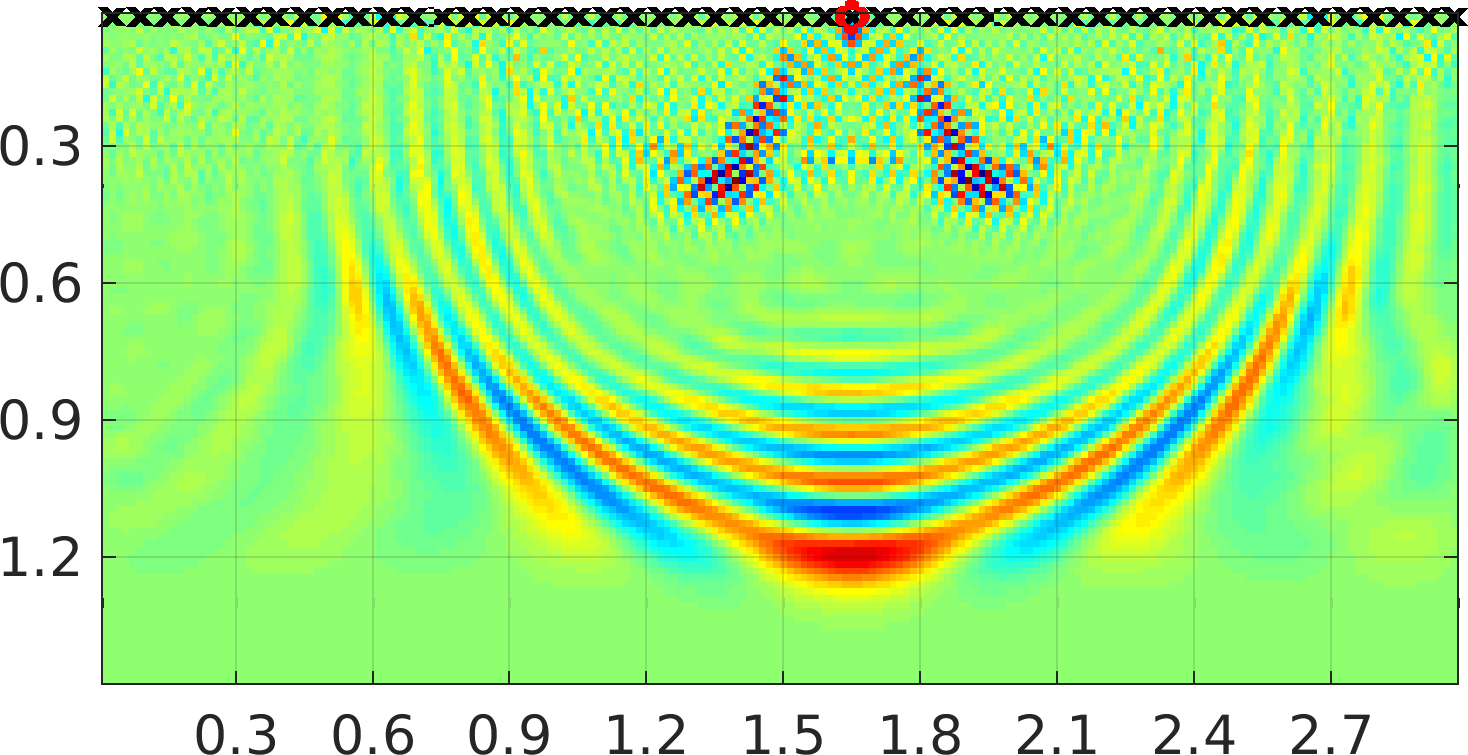} &
\includegraphics[width=0.22\textwidth]{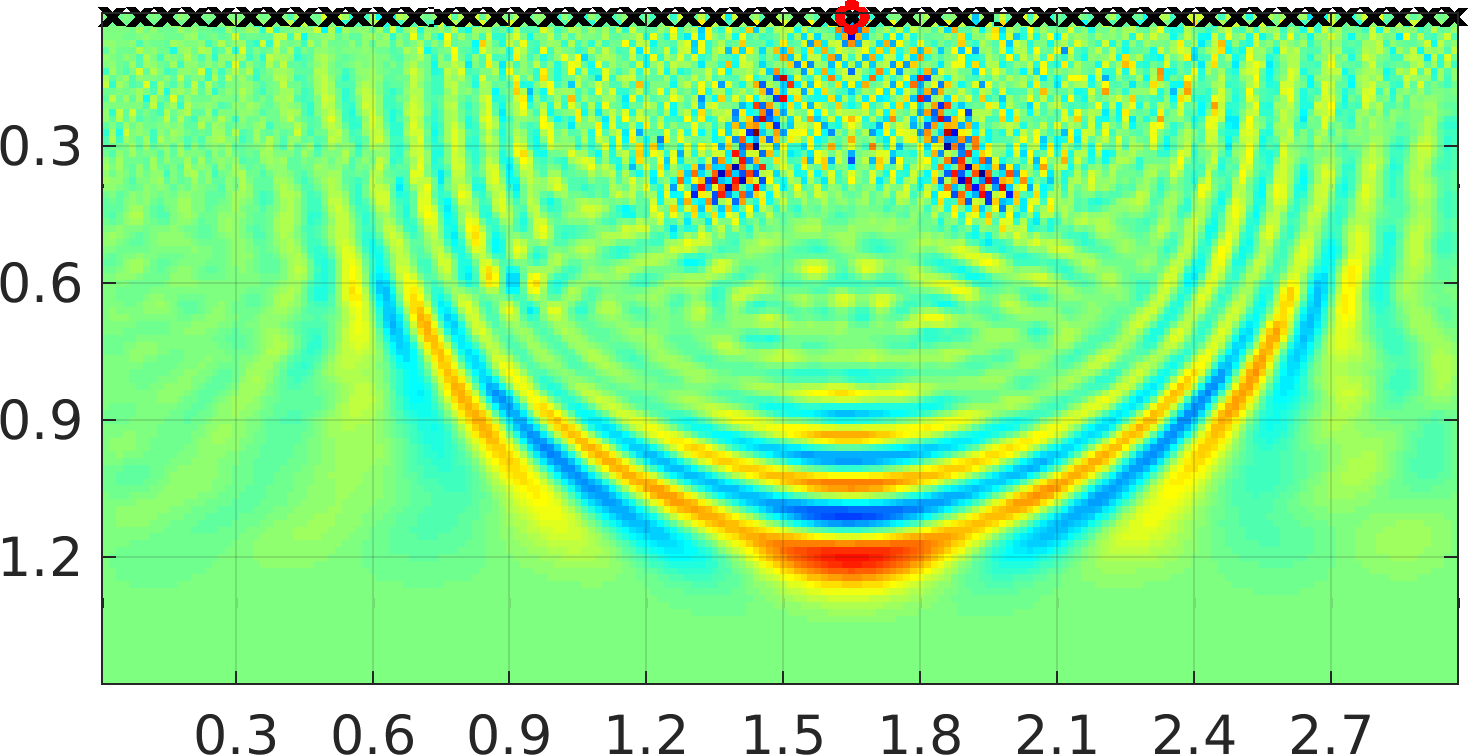}
\end{tabular}
\vspace{-0.12in}\caption{Primary snapshots and their orthonormal
  counterparts for the reference medium with $\sigma^0(\bx) \equiv 1$
  and the true medium.  For every snapshot we plot a single column of
  the $N \times mn$ matrices, corresponding to the source
  $j=28$ (red $\circ$) out of $m=50$ sensors in the array (black
  $\times$). The rows (top to bottom) correspond to times
  $t_k = k \tau$, for $k = 25, 30$. The axes are in km units.}
\label{fig:UV}
\end{figure}

In Figure~\ref{fig:UV} we plot the primary snapshots at two time
instances $t_k = k \tau$, with $k = 25$ and $k = 30$, for the
reference medium with $\sigma^0 \equiv 1$ and the scattering medium
displayed in Figure \ref{fig:model2d}. We also display the orthonormal
snapshots.  The primary snapshots in the reference medium (first
column) show the wavefront.  In the scattering medium (second column)
there are multiple reflections behind the wavefront.  The orthonormal
snapshots for both media (third and forth columns) have a
``smile-like'' shape with a ``thick lower lip''. They can be viewed as
approximations of delta functions. The reflections are suppressed
in the orthonormal snapshots, and we note that they are almost the
same in both media i.e., they are almost independent of $\sigma$.  A
similar result holds for the dual snapshots, not shown here.

\begin{figure}
\begin{tabular}{ccc}
Scattered data & Born approximation & DtB  \\
\includegraphics[width=0.3\textwidth]{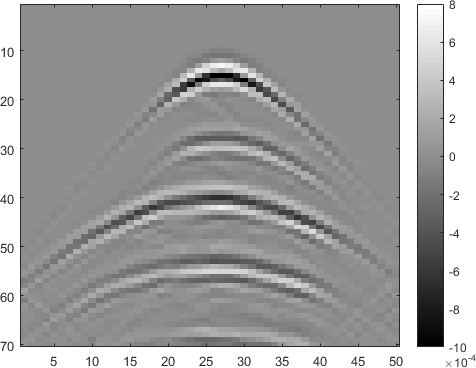} & 
\includegraphics[width=0.3\textwidth]{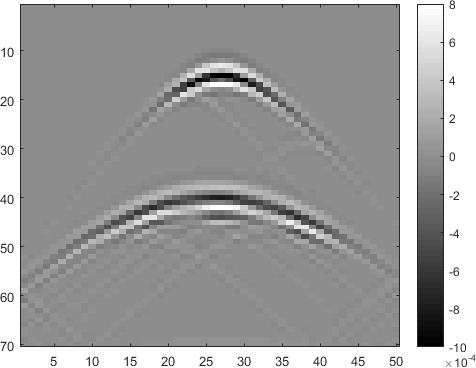} &
\includegraphics[width=0.3\textwidth]{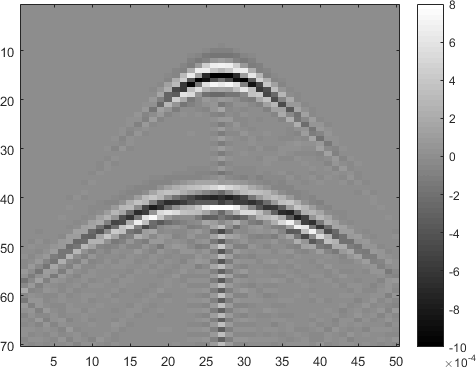}   
\end{tabular}
\vspace{-0.12in}\caption{From left to right: raw scattered data, Born approximation
  and the data transformed by the DtB algorithm. The abscissa is the
  sensor index, and the ordinate is the index $k$ of
  the discrete time instants $t_k = \tau k$, $\tau=0.034s$.}
\label{fig:DtB_1st_2d}
\end{figure}

\begin{figure}
\begin{center}
\includegraphics[width=0.75\textwidth]{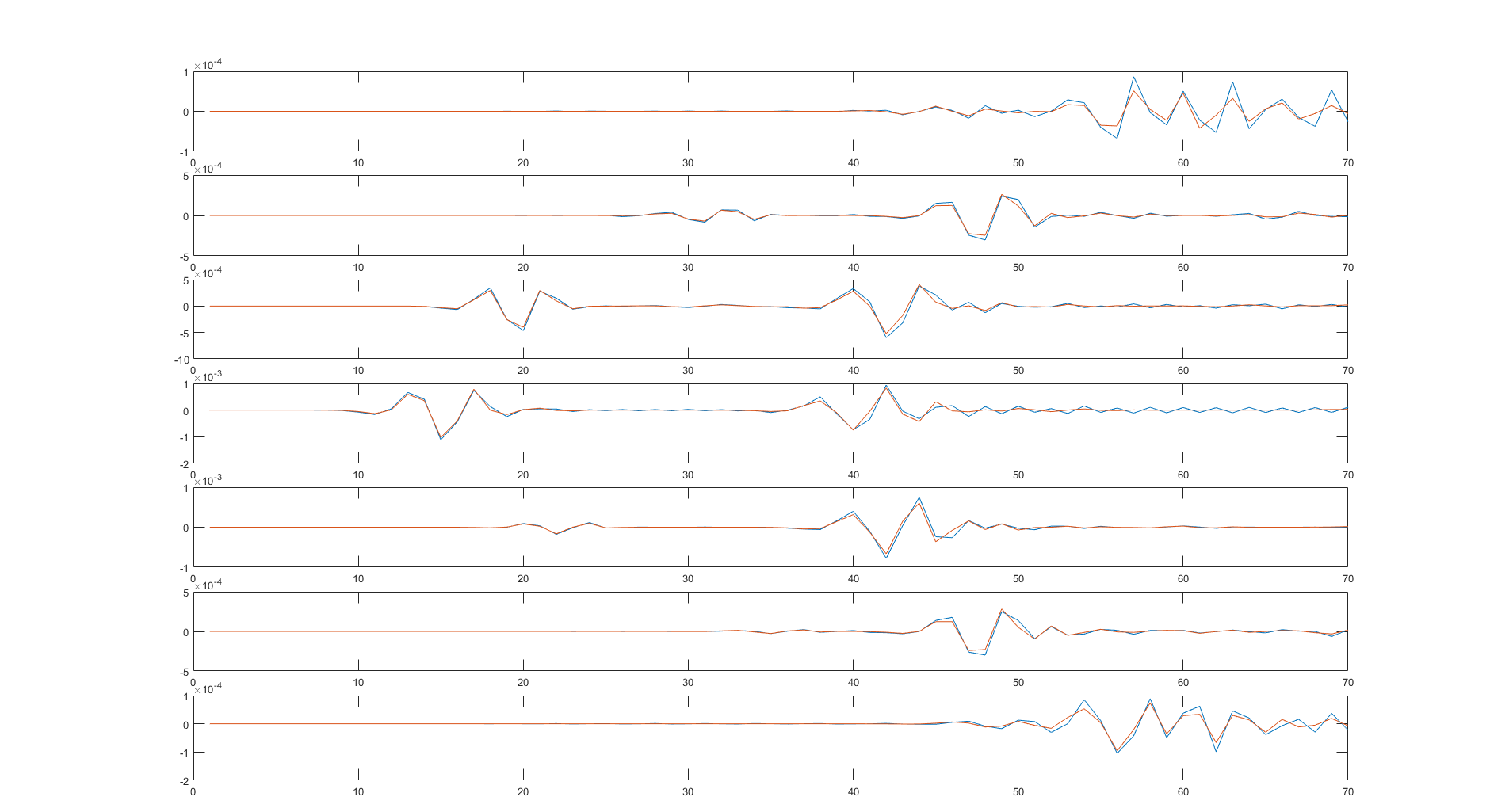} 
\vspace{-0.12in}\caption{Comparison of a few traces of the Born approximation (red)
  and the DtB algorithm output (blue) from
  Figure~\ref{fig:DtB_1st_2d}. The abscissa is the time index.}
\label{fig:DtBvsBorn}
\end{center}
\end{figure}

In the left plot of Figure~\ref{fig:DtB_1st_2d} we display the raw
scattered data at the sensors, due to the excitation from the sensor
at the center of the linear array, lying just below the top boundary.
The Born approximation and the data transformed by the DtB algorithm
are in the middle and right plots. The results are almost the same. To
illustrate better the agreement between the Born approximation and the
output of the DtB algorithm, we display in Figure~\ref{fig:DtBvsBorn}
a comparison of several traces (signals at certain receivers) from
Figure~\ref{fig:DtB_1st_2d}.

\begin{figure}
\begin{tabular}{cc}
Image with raw data & Image with DtB transformed data \\
\includegraphics[width=0.45\textwidth]{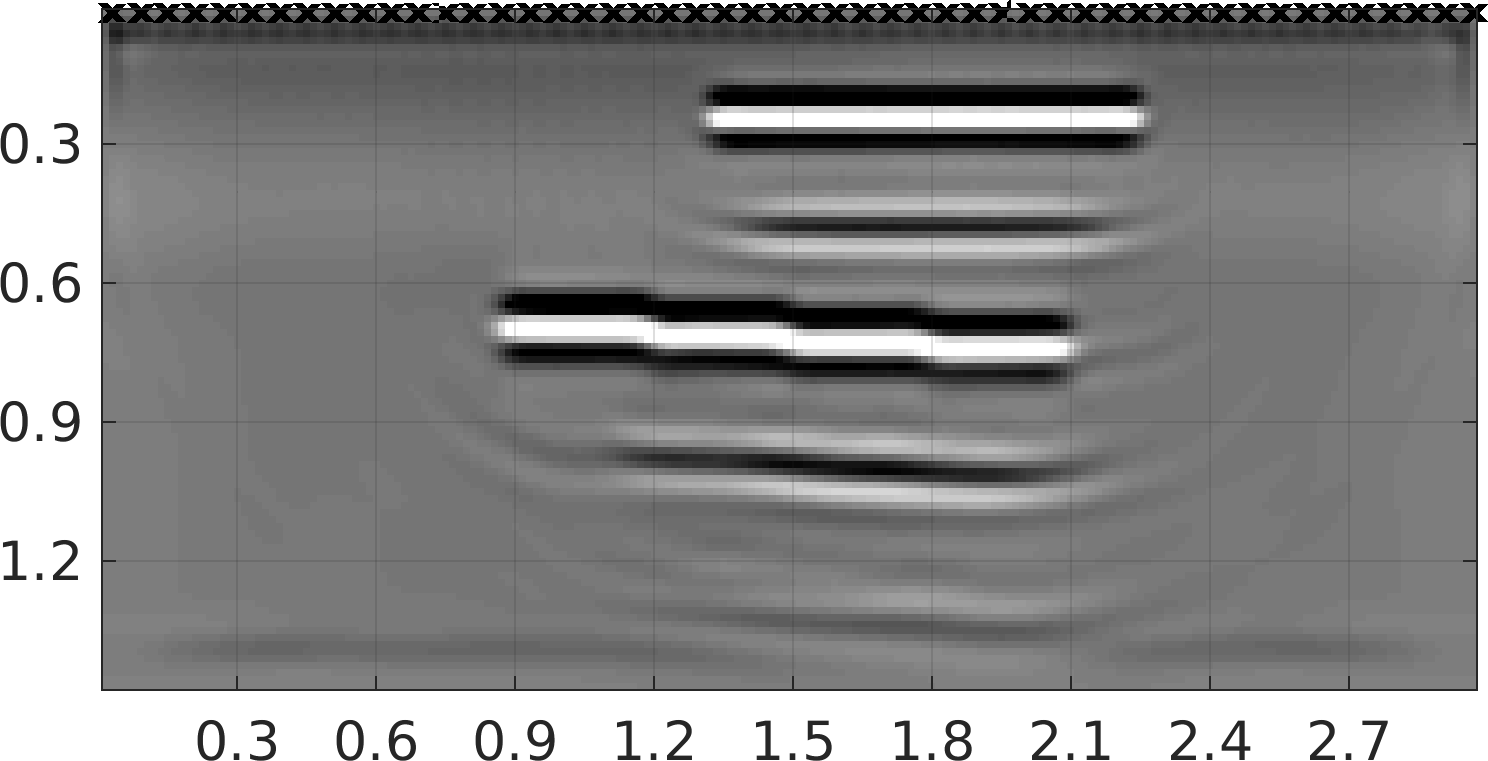} & 
\includegraphics[width=0.45\textwidth]{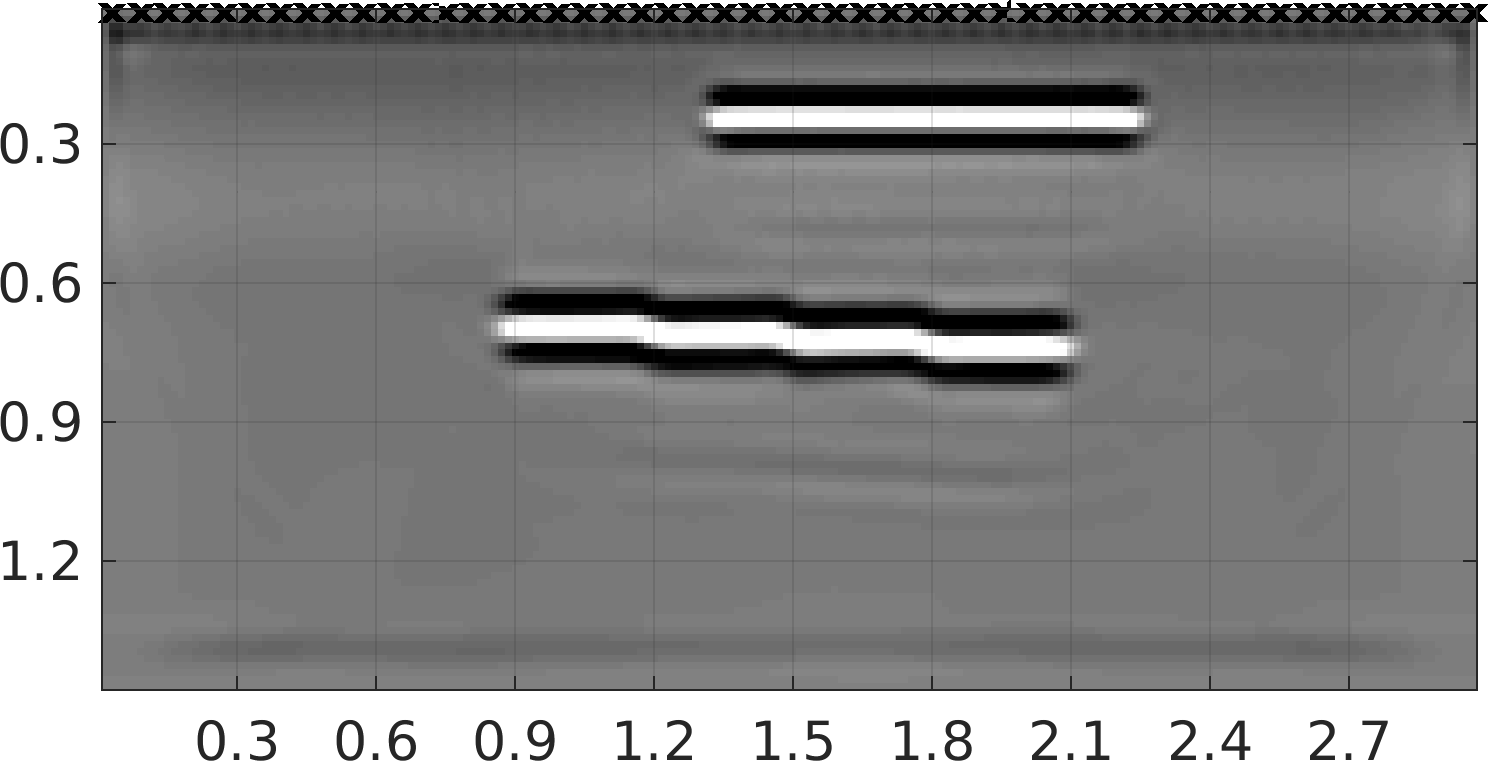} \end{tabular}
\vspace{-0.12in}\caption{Reverse time migration image obtained with the raw data
(left) and the data transformed with the DtB algorithm (right). The axes are in km units, 
as in Figure \ref{fig:model2d}.}
\label{fig:Images}
\end{figure}

To illustrate the benefit of the DtB transformation  on imaging, 
we display in Figure \ref{fig:Images} the reverse time migration images\footnote{We refer to \cite{Biondi,symes2009seismic,symes1995mathematics}
for details on the reverse time migration. It amounts to taking the data, time reversing it 
and backpropagating it in the reference (nonscattering) medium, with velocity 
$c(\bx)$. The image displays the resulting wave field evaluated at the travel time to points in the imaging region. } 
obtained with the raw data and the transformed data shown in the left and right plots of 
Figure \ref{fig:DtB_1st_2d}.  The artifacts due to multiple scattering are evident in the image shown on the 
left, which displays multiple ghost reflectors. The image obtained with the transformed data,
shown on the right, does not have multiple artifacts and localizes well the two reflectors.

\begin{figure}
\begin{tabular}{ccc}
Model $c(\bx)$ & Scattered data  & DtB  \\
\includegraphics[width=0.33\textwidth, height=3cm]{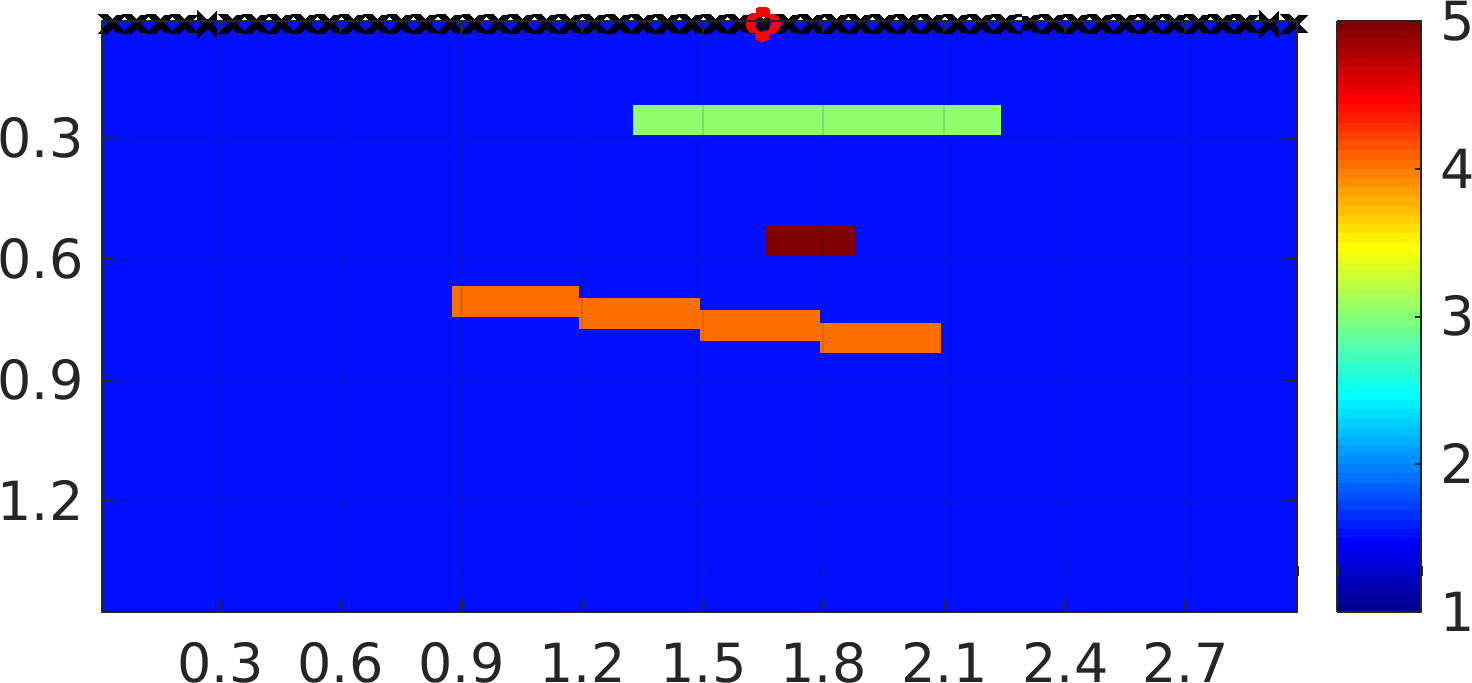} & 
\includegraphics[width=0.3\textwidth]{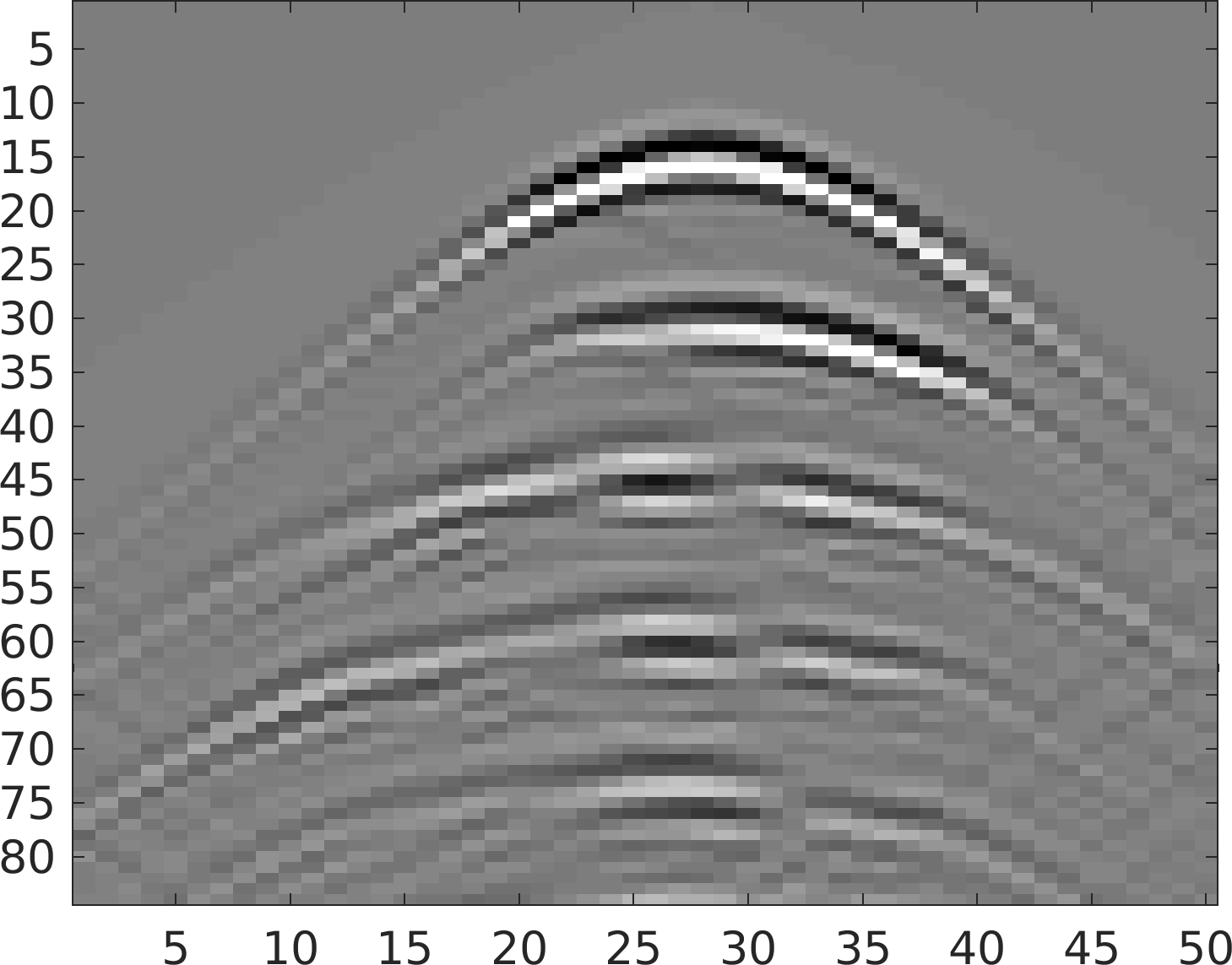} &
\includegraphics[width=0.3\textwidth]{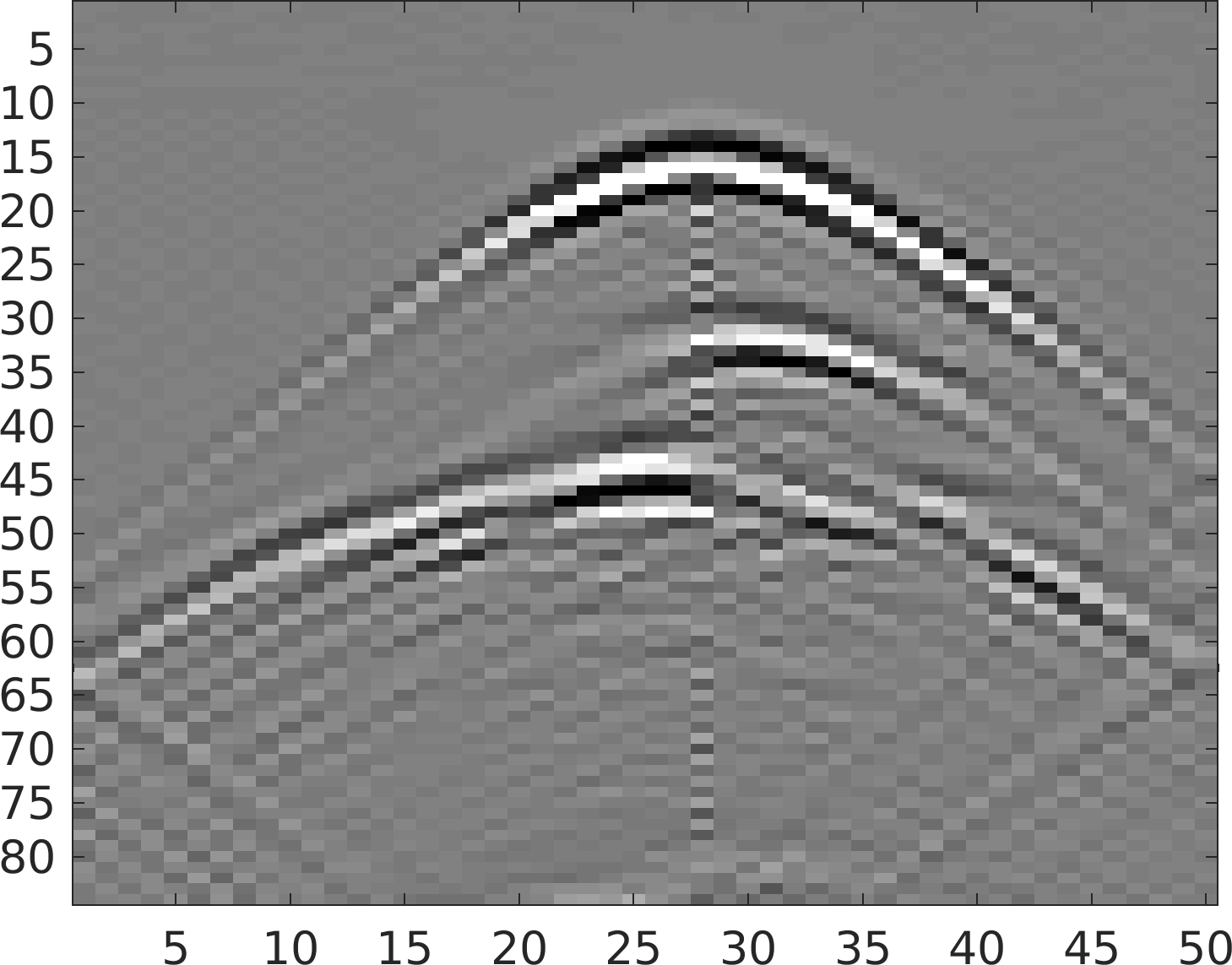}   
\end{tabular}
\vspace{-0.11in}
\caption{ Left: model of a medium with three inclusions.  The value of
  the wave speed $c(\bx)$ is shown in the color bar, in units of
  km/s. The axes are in km units and the array sits on the top, as in the previous simulations. Middle: Data gathered for a single source
  shown as red $\circ$ in the left figure.   The abscissa is the sensor receiver index and the ordinate
  is the discrete time sample with step $\tau=0.0225s$.  Right: the result of the 
  DtB algorithm. }
\label{fig:motive}
\end{figure}

As we mentioned in the introduction, in practice only the smooth part of 
$c(\bx)$ may be known. To illustrate that the DtB algorithm can 
deal with perturbations of the sound speed, we present in 
Figure \ref{fig:motive} numerical results for  
three inclusions embedded in a medium with constant wave speed $c_0 = 1$km/s and constant density 
$\rho$. The inclusions are modeled by the variation of $c(\bx)$ displayed in the 
left plot, but the density is kept constant i.e., $\sigma(\bx) = \rho c(\bx)$. Only $c_0$ is assumed known
in the DtB algorithm, meaning that  we used the incorrect speed $c_0$ instead of the 
true $c(\bx)$. The data gathered by the array, for the excitation from the source shown with a red circle in the left plot,  
are displayed in the middle plot. They contain the primary reflections from each inclusion 
and multiply scattered reflections between the inclusions. The output of the 
DtB algorithm is displayed in the right plot of Figure \ref{fig:motive}. The 
multiply scattered echoes are removed and there are three, clearly separated echoes, 
corresponding to each inclusion. Note  the unmasking of the second echo, 
due to the smaller inclusion, that was mixed with a multiply scattered echo in the middle 
plot.

\section{Summary}
\label{sect:sum}
This paper is motivated by the inverse scattering problem for the wave
equation, where an array of sensors probes an unknown scattering
medium with pulses and measures the reflected waves.  The goal of the
inversion is to estimate the perturbations of the acoustic impedance
in the medium, which cause wave scattering. We introduced a direct,
linear-algebra based algorithm, called the Data to Born (DtB)
algorithm, for transforming the data collected by the array to data
corresponding to the single scattering (Born) approximation. 
These data can then be used by any  off-the-shelf algorithms that  incorporate state of the art 
linear inversion methodologies.  The key
ingredient in the DtB algorithm is a data driven, reduced order model
(ROM), that approximates the wave propagator operator. Because
the DtB algorithm involves only linear algebra operations, like
matrix-matrix multiplications and block Cholesky factorizations, the
cost of the algorithm is $O((mn)^3)$, where $m$ is the number of
sensors and $2n$ is the number of time samples in the measurements.
{ However, due to the Toeplitz-plus-Hankel structure of the mass and stiffness matrices,  this cost can be reduced, possibly to 
$O(m^3n^2\log n)$. This will be a subject of future research.}
\section*{Acknowledgments}
\label{sec:acknowledge}
This material is based upon research supported in part by the
U.S. Office of Naval Research under award number N00014-17-1-2057 to
Borcea and Mamonov. Mamonov was also partially supported by
the NSF grant DMS-1619821.

\appendix
\section{The tridiagonal structure of the ROM propagator}
\label{ap:A}
In this appendix we show that the the ROM propagator given by the
projection \eqref{eq:1D57} is a tridiagonal matrix. Obviously,
$\tilde{\bbP}$ is symmetric, so it suffices to show that its entries
\begin{equation}
\label{eq:A1}
(\tilde{\bbP})_{j+l,j} = \bV_{j+l}^T \bbP \bv_j
\end{equation}
are zero when $l \ge 2$. 

We obtain from equation \eqref{eq:1D55p} that $ \displaystyle \bv_j =
\sum_{i=0}^{j-1} (\bR^{-1})_{i+1,j} \bP_i, $ where we used that the
inverse $\bR^{-1}$ of the upper triangular matrix $\bR$ is upper
triangular.  The relation \eqref{eq:ChebT4} satisfied by the Chebyshev
polynomials and definition \eqref{eq:solSnaps} give
\begin{equation}
\bbP \bP_i = \mathcal{T}_1(\bbP) \mathcal{T}_i(\bbP) \bb = \frac{1}{2}
\Big[\mathcal{T}_{i+1}(\bbP) + \mathcal{T}_{|i-1|}(\bbP) \Big] \bb =
\frac{1}{2} \big(\bP_{i+1} + \bP_{|i-1|}\big),
\label{eq:A3}
\end{equation}
so equation \eqref{eq:A1} becomes 
\begin{equation}
(\tilde{\bbP})_{j+l,j} = \frac{1}{2} \sum_{i=0}^{j-1} (\bR^{-1})_{i+1,j}  \Big( 
\bV_{j+l}^T \bP_{i+1} + \bV_{j+l}^T \bP_{|i-1|}\Big).
\label{eq:A4}
\end{equation}
Each term in this sum can be calculated from \eqref{eq:1D55} as $
\bv_{j+l}^T \bP_i = (\bR)_{j+l,i+1}, $ so we obtain
\begin{equation}
(\tilde{\bbP})_{j+l,j} = \frac{1}{2} \sum_{i=0}^{j-1} (\bR^{-1})_{i+1,j}  \Big[ 
(\bR)_{j+l,i+2} + (\bR)_{j+l,|i-1|+1}\Big).
\label{eq:A6}
\end{equation}
Note that $j+l > \max\{i+2,|i-1|+1\}$ when $l \ge 2$ and $i = 0,
\ldots, j-1$, so the right hand side in \eqref{eq:A6} is zero by the
upper triangular structure of $\bR$. This means that $\tilde{\bbP}$ is
tridiagonal.

\section{Computation of the block-bidiagonal   $\tbL_q$} 
\label{ap:blockCholesky}

{We describe the computation of the block-Cholesky factor
  $\tbL_q$ using an approach outlined in \cite{DMZmultiscale2017}.}
As mentioned in section \ref{sect:mimoromcalc}, the block Cholesky
factorization is not uniquely defined.  Clearly, if $\bQ_k={\itbf
  I}_m,$ the diagonal blocks $\bR_{k,k},$ are symmetric, for $k =
1, \ldots, n.$ We denote the corresponding MIMO ROM matrix given by
\eqref{eq:1D56p} as $\tilde{\bbP}^{{\itbf I}}$.

For non-trivial orthogonal matrices $\bQ_k$, the MIMO ROM matrix
$\tilde{\bbP}^{\bQ}$ is given by $\bQ\tilde{\bbP}^{{\itbf I}}\bQ^T$
where $\bQ$ is the block-diagonal matrix with orthogonal blocks
$\bQ_k$, $k=1,\ldots, n$. However, the block bidiagonal factor
$\tbL_q$ in \eqref{eq:Chol} has to be consistent with the matrix
analogue of the recursion \eqref{eq:O1}--\eqref{eq:O3p} for computing
the primary and dual orthogonalized block snapshots $\overline{\bf
  p}_{j}\in \mathbb{R}^{N \times m}$ and $\overline{\bf u}_{j}\in
\mathbb{R}^{N \times m}$,
\begin{align}
\bgamma^{-1}_j\left(\overline{\bf p}_{j+1} - \overline{\bf p}_{j}\right) = - \cbL_q
\overline{\bf u}_j, \quad 
\bgammahat^{-1}_j\left(\overline{\bf u}_{j} -
\overline{\bf u}_{j-1}\right) = \cbL_q^T \overline{\bf p}_j, \quad j \ge 1, 
\label{eq:O1b}
\end{align}
with initial conditions $ \overline{\bf p}_{1} = \bb$
and $\overline{\bf u}_{0} = 0,$ and symmetric positive definite matrix
coefficients
\begin{equation}
\hspace{-0.1in} 0 < \bgamma_j = \bgamma^T_j = \Big(
\overline{\bf u}_j^T \overline{\bf u}_j \Big)^{-1} \in \mathbb{R}^{m
  \times m}, \quad 0 < \bgammahat_j = \bgammahat^T_j = \Big(
\overline{\bf p}_j^T \overline{\bf p}_j \Big)^{-1} \in \mathbb{R}^{m
  \times m},
\label{eq:O3b}
\end{equation}
for $j = 1, \ldots, n$. Then,
\begin{align}
(\tbL_q)_{j,j} &= - \sqrt{\bgammahat^{-1}_j}\sqrt{\bgamma^{-1}_j }, ~ ~ 1 \le j \le n, 
  \quad (\tbL_q)_{j+1,j} =
  \sqrt{\bgammahat^{-1}_{j+1}}\sqrt{\bgamma^{-1}_{j}} ~ ~ 1 \le j \le n-1.
\label{eq:O3pb}
\end{align}
We now determine the matrix $\bQ$ such that the factorization 
\begin{equation}
\label{eq:romeqchol}
\xi \Big( \bQ \tilde{\bbP}^{{\itbf I}}\bQ^T \Big) = 
\bQ \xi \Big( \tilde{\bbP}^{{\itbf I}} \Big) \bQ^T = 
\tbL_q\tbL^T_q,
\end{equation}
corresponds to $\tbL_q$ of the form \eqref{eq:O3pb} with symmetric
positive definite $\bgamma_j$ and $\bgammahat_j$, $j=1,\ldots, n$.

Denote the diagonal and the off-diagonal blocks of $\xi \Big(
\tilde{\bbP}^{{\itbf I}} \Big)$ by $\balpha_j \in \mathbb{R}^{m \times
  m}$, for $j = 1,\ldots,n$ and $\bbeta_j \in \mathbb{R}^{m \times
  m}$, for $j = 2,\ldots,n$. By definition,
\begin{equation}
\label{eq:gammahat1}
\bgammahat_1= (\bb^T \bb)^{-1}.
\end{equation}

The remaining matrix coefficients $\bgamma_j$ and $\bgammahat_j$ are
obtained from \eqref{eq:romeqchol} block-wise. From the first diagonal
block we obtain that $\sqrt{\bgammahat^{-1}_1} \bgamma^{-1}_1
\sqrt{\bgammahat^{-1}_1} = \bQ_{1} \balpha_1 \bQ^T_1$ or,
equivalently,
\begin{equation}
\label{eq:gamma1}
\bgamma_1 = \left(\sqrt{\bgammahat_1}\bQ_1 \balpha_1 \bQ^T_1
\sqrt{\bgammahat_1}\right)^{-1}.
\end{equation}

Clearly, $\bgamma_1=\bgamma^T_1$ for any matrix $\bQ_1$, so for
simplicity we set $\bQ_1={\itbf I}_m$.  Then, from the off-diagonal
blocks for $1 \le j\le n-1$ we have $\sqrt{\bgammahat^{-1}_j}
\bgamma^{-1}_j \sqrt{\bgammahat^{-1}_{j+1}} = \bQ_{j}
\bbeta_{j+1}\bQ^T_{j+1}$.  Hence,
$\sqrt{\bgammahat^{-1}_{j+1}}\bQ_{j+1} = \bgamma_j
\sqrt{\bgammahat_j}\bQ_{j} \bbeta_{j+1}$. That is to say, the pair of
matrices $\sqrt{\bgammahat^{-1}_{j+1}}$ and $\bQ_{j+1}$ is a (left)
polar decomposition of $ \bM_j = \bgamma_j \sqrt{\bgammahat_j} \bQ_j
\bbeta_{j+1}.  $ Its solution is
\begin{align}
\bgammahat_{j+1} & = \left( \bM_j \bM^T_j \right)^{-1}, \quad
\bQ_{j+1} = \sqrt{\bgammahat_{j+1}} \bM_j.
\label{eq:poldec}
\end{align}
Finally, considering the diagonal blocks for $1 \le j \le n-1$, we obtain 
\[
\sqrt{\bgammahat^{-1}_{j+1}} \left( \bgamma^{-1}_j + \gamma^{-1}_{j+1} \right) 
\sqrt{\bgammahat^{-1}_{j+1}} = \bQ_{j+1} \balpha_{j+1} \bQ^T_{j+1}
\] 
and therefore 
\begin{equation}
\label{eq:gammaip1}
\bgamma_{j+1} = \left( \sqrt{\bgammahat_{j+1}} \bQ_{j+1} \balpha_{j+1}
\bQ^T_{j+1} \sqrt{\bgammahat_{j+1}} - \bgamma^{-1}_j \right)^{-1}.
\end{equation}

\vspace{0.05in}\begin{alg}[Computation of $\tbL_q$]
\label{alg:blockholwgm}~\\
\textbf{Input:} the block tridiagonal matrix $\tilde{\bbP}^{{\itbf I}}
\in \mathbb{R}^{nm \times nm}$ with $m \times m$ blocks and $\bb^T
\bb$.

\vspace{0.05in} \noindent To find a block diagonal $\bQ$ such that the
block Cholesky factorization $\tbL_q\tbL^T_q$ of {$\xi \Big(
  \tilde{\bbP}^{\bQ} \Big)$} has factors $\tbL_q$ in the form
\eqref{eq:O3pb}, perform the following steps:

\vspace{0.05in} \noindent 1.  Compute $\bgammahat_1$ via
\eqref{eq:gammahat1} and $\bgamma_1$ via \eqref{eq:gamma1} for
arbitrary $\bQ_1$ (for simplicity, we set $\bQ_1={\itbf I}_m$)

\vspace{0.05in} \noindent 2. For $j = 1,\ldots,n-1:$ Compute
$\bgammahat_{j+1}$ and $\bQ_{j+1}$ via \eqref{eq:poldec} and
$\bgamma_{j+1}$ via \eqref{eq:gammaip1}.

\vspace{0.05in} \noindent
3.  Compute $\tbL_q$ via \eqref{eq:O3pb}

\vspace{0.05in} \noindent
\textbf{Output:} the block diagonal orthogonal matrix $\bQ$, the block
tridiagonal propagator matrix $\tilde{\bbP}^{\bQ}=\bQ\tilde{\bbP}^{{\itbf
    I}}\bQ^T$ and the block lower bidiagonal factor $\tbL_q$ of $\xi
\Big( \tilde{\bbP}^{\bQ} \Big)$ consistent with \eqref{eq:O3pb}.
\end{alg}

\bibliography{biblio} \bibliographystyle{siam}
\end{document}